\documentclass[preprint,12pt]{elsarticle}

\usepackage{graphicx}
\usepackage{amsmath}
\usepackage{amssymb}
\usepackage{amsthm}

\usepackage{algorithm}
\usepackage{algpseudocode}
\algrenewcommand\textproc{}

\usepackage{caption}
\usepackage{xcolor}

\usepackage{subfig}

\usepackage{ulem}

\usepackage{multicol}
\usepackage{multirow}
\usepackage{adjustbox}

\journal{Elsevier}

\usepackage[margin=1in]{geometry}

\usepackage{graphicx}
\usepackage{pstool}
\usepackage{wrapfig}
\usepackage{caption}
\usepackage[section]{placeins} % Defines \FloatBarrier which does not allow
\usepackage{fancyhdr} % Required for custom headers
\usepackage{lastpage} % Required to determine the last page for the footer
\usepackage{extramarks} % Required for headers and footers

\usepackage{xcolor} % Colors
\usepackage{enumerate} % Additional options for \enumerate{}
\usepackage{paralist} % Additional options for lists; inline lists
\usepackage{amsmath, amsthm, amssymb, mathtools}
\usepackage{mathabx, pifont, stmaryrd} % More math symbols, fonts
\usepackage[explicit]{titlesec} % Control over section/subsection styles
\usepackage{etoolbox} % Boolean variables and more
\usepackage{bibentry} % Bibliography entries anywhere in doc
\makeatletter\let\saved@bibitem\@bibitem\makeatother % Do not remove if bibentry
\usepackage[colorlinks, bookmarksopen, bookmarksnumbered,
            citecolor=red,urlcolor=red]{hyperref} % Hyperlinks
\makeatletter\let\@bibitem\saved@bibitem\makeatother % Do not remove if bibentry
\usepackage{algpseudocode}

\algnewcommand{\algorithmicand}{\textbf{ and }}
\algnewcommand{\algorithmicor}{\textbf{ or }}
\algnewcommand{\OR}{\algorithmicor}
\algnewcommand{\AND}{\algorithmicand}
\algnewcommand{\var}{\texttt}

\newtheorem{remark}{Remark}

\usepackage{bm} % required for: \bm{}
\usepackage{amsfonts} % required for: \mathbb{}, \mathcal{}, ...

\newcommand{\Ical}{\ensuremath{\mathcal{I}}}

\newcommand{\Nbb}{\ensuremath{\mathbb{N} }}

\newcommand{\Rbb}{\ensuremath{\mathbb{R} }}

\usepackage{tikz}
\usepackage{pgfplots}
\usepackage{pgfplotstable, filecontents, booktabs}
\pgfplotsset{compat=1.9}

\usetikzlibrary{pgfplots.groupplots}
\usepgfplotslibrary{fillbetween}
\usetikzlibrary{calc,fit,matrix,arrows,automata,positioning,shapes}
\usetikzlibrary{arrows.meta}

\pgfplotsset{select coords between index/.style 2 args={
    x filter/.code={
        \ifnum\coordindex<#1\fi
        \ifnum\coordindex>#2\fi
    }
}}

\tikzset{
 invisible/.style={opacity=0},
 visible on/.style={alt={#1{}{invisible}}},
 alt/.code args={<#1>#2#3}{%
   \alt<#1>{\pgfkeysalso{#2}}{\pgfkeysalso{#3}}
 },
}

\begin{document}

\begin{frontmatter}

%% Title, authors and addresses

%% use the tnoteref command within \title for footnotes;
%% use the tnotetext command for theassociated footnote;
%% use the fnref command within \author or \address for footnotes;
%% use the fntext command for theassociated footnote;
%% use the corref command within \author for corresponding author footnotes;
%% use the cortext command for theassociated footnote;
%% use the ead command for the email address,
%% and the form \ead[url] for the home page:
%% \title{Title\tnoteref{label1}}
%% \tnotetext[label1]{}
%% \author{Name\corref{cor1}\fnref{label2}}
%% \ead{email address}
%% \ead[url]{home page}
%% \fntext[label2]{}
%% \cortext[cor1]{}
%% \affiliation{organization={},
%%             addressline={},
%%             city={},
%%             postcode={},
%%             state={},
%%             country={}}
%% \fntext[label3]{}

\title{An adaptive, training-free reduced-order model for convection-dominated problems based on hybrid snapshots}
%\title{An adaptive reduced-order model for convection-dominated problems based on hybrid snapshot computation}

%% use optional labels to link authors explicitly to addresses:
%% \author[label1,label2]{}
%% \affiliation[label1]{organization={},
%%             addressline={},
%%             city={},
%%             postcode={},
%%             state={},
%%             country={}}
%%
%% \affiliation[label2]{organization={},
%%             addressline={},
%%             city={},
%%             postcode={},
%%             state={},
%%             country={}}

\author[inst1]{Victor Zucatti\fnref{fn1}}
\ead{vzucatti@nd.edu}

\author[inst1]{Matthew J. Zahr\fnref{fn2}\corref{cor1}}
\ead{mzahr@nd.edu}

\affiliation[inst1]{organization={University of Notre Dame},
            city={Notre Dame},
            postcode={IN 46556},
            country={United States of America}}
\cortext[cor1]{Corresponding author}

\fntext[fn1]{Graduate Student, Department of Aerospace and Mechanical
             Engineering, University of Notre Dame}
\fntext[fn2]{Assistant Professor, Department of Aerospace and Mechanical
             Engineering, University of Notre Dame}

\begin{abstract}
The vast majority of reduced-order models (ROMs) first obtain a low dimensional representation of the problem from high-dimensional model (HDM) training data which is afterwards used to obtain a system of reduced complexity. Unfortunately, convection-dominated problems generally have a slowly decaying Kolmogorov $n$-width, which makes obtaining an accurate ROM built solely from training data very challenging. The accuracy of a ROM can be improved through enrichment with HDM solutions; however, due to the large computational expense of HDM evaluations for complex problems, they can only be used parsimoniously to obtain relevant computational savings. In this work, we exploit the local spatial coherence often exhibited by these problems to derive an accurate, cost-efficient approach that repeatedly combines HDM and ROM evaluations without a separate training phase. Our approach obtains solutions at a given time step by either fully solving the HDM or by combining partial HDM and ROM solves. A dynamic sampling procedure identifies regions that require the HDM solution for global accuracy and the reminder of the flow is reconstructed using the ROM. Moreover, solutions combining both HDM and ROM solves use spatial filtering to eliminate potential spurious oscillations that may develop. We test the proposed method on inviscid compressible flow problems and demonstrate speedups up to a factor of five.
\end{abstract}

%%Graphical abstract
% \begin{graphicalabstract}
% \includegraphics{grabs}
% \end{graphicalabstract}

%%Research highlights
%\begin{highlights}
%\item Research highlight 1
%\item Research highlight 2
%\end{highlights}

\begin{keyword}
%% keywords here, in the form: keyword \sep keyword
adaptive model reduction \sep proper orthogonal decomposition \sep hyperreduction \sep sparse sampling \sep convection-dominated problems
%% PACS codes here, in the form: \PACS code \sep code
% \PACS 0000 \sep 1111
%% MSC codes here, in the form: \MSC code \sep code
%% or \MSC[2008] code \sep code (2000 is the default)
% \MSC 0000 \sep 1111
\end{keyword}

\end{frontmatter}

% \section*{Possible titles}

% \begin{itemize}
%     \item POD-Based Hybrid Snapshot Adaptive Model Reduction to Accelerate Unsteady Problems
%     % 
%     \item Hybrid snapshot adaptive model reduction to accelerate large-scale unsteady problems
%     % 
%     \item Accelerating large-scale simulations via hybrid snapshot computations: adaptive model reduction for convection-dominated problems
%     % 
%     \item An adaptive reduced-order model for large-scale convection-dominated problems based on hybrid snapshot computation
% \end{itemize}

\section{Introduction}
\label{sec:intro}

Today's computational power enables the numerical solution of complex engineering problems; however, these computations can easily require hundreds of millions of degrees of freedom to produce accurate results \cite{RICCIARDI2021115933} and, thus, high-fidelity many-query analyses are still impractical in many scenarios such as design optimization, flow control and uncertainty quantification, to name a few.
Fortunately, the large amount of data generated by high-dimensional models (HDMs) can be used to build a reduced-order model (ROM).
A two-step (offline/online) approach is the standard when building ROMs for time-dependent problems. In the offline stage, a smaller dimensional representation is obtained from HDM training data and used to generate a lower complexity model through physics-based \cite{Rowley2004,Kevin01,Kevin02_comp} or data-driven methods \cite{2016_brunton_sindy,2016_peherstorfer_opInf,lui_wolf_2019}. This is a precomputation step performed only once but can be very costly given the high-dimensional data dependency. On the other hand, the online stage consists of solving the resulting system of equations of reduced dimensionality (e.g., up to four orders of magnitude smaller \cite{lui_wolf_2019}).
Unfortunately, despite the considerable research done in the last 20 years, ROMs still suffer from a multitude of problems (e.g., instability, inaccuracy, failure to generalize beyond training) making them generally unreliable in an industrial setting \cite{Rowley2004,Cazemier1998,Noack_01_pressure}.
This is particularly the case when modeling time-dependent convection-dominated problems such as those usually found in viscous or high-speed computational fluid dynamics (CFD) problems.
Multiple correction methods have been proposed \cite{bergmann_enablers,iliescu_01_burgers,iliescu_02_turbulent,KALASHNIKOVA_2014,grimberg2020stability,zucatti_cmame} and have been rather successful in improving ROM stability. However, they have done very little to improve ROM predictive capabilities.

For convection-dominated problems, failure to generalize has been mainly attributed to the slowly decaying Kolmogorov $n$-width of linear subspace approximations \cite{2016_Ohlberger}. This is also sometimes referred to as Kolmogorov barrier because the error slowly decaying with the dimension of the reduced space limits the achievable accuracy of ROMs in practice and requires a substantial amount of training data, which can be infeasible to collect offline. The Kolmogorov barrier can be overcome, for example, by the use of nonlinear model reduction techniques.
In \cite{kevin_04_autoencoder}, a nonlinear manifold is obtained through deep convolutional autoencoders and combined with projection-based methods to produce ROMs capable of outperforming their linear counterparts. 
Quadratic manifolds have been used with both physics-based \cite{2022_Barnett} and data-driven \cite{2022_Geelen_quad} methods for order reduction.
Alternatively, nonlinear manifolds have been constructed by composing a traditional subspace approximation with a transformation to the underlying domain, which has proven particularly effective for shock-dominated problems \cite{2020_Tommaso,2021_Mirhoseini}.
Another solution is to exploit the local low-rank structure of this class of problems \cite{2020_Peherstorfer_AADEIM}. In \cite{Amsallem2012_local,Amsallem2015_local}, local low-rank subspaces are systematically obtained by partitioning of the state space. Results show that local subspaces improves ROMs accuracy and speed by reducing the dimensionality of each subspace.

Adaptive reduced-order models (AROMs) \cite{Bai2020, Bai2022,Feng2021,2015_Peherstorfer_ADEIM,2020_Peherstorfer_AADEIM,2022_uy_aadeim,huang2023predictive} provide a different approach by continuously combining HDM and ROM operations. 
Predictive capabilities can be improved by alternating between HDM and ROM generated snapshots \cite{Bai2020, Bai2022,Feng2021}. 
In \cite{Bai2020, Bai2022}, on-the-fly criteria relying on the reduced basis sufficiency is used to determine when to use the HDM or local ROM. If deemed necessary, fast low-rank singular value decomposition (SVD) modifications \cite{BRAND200620} are used to update the reduced-order basis. This methodology was successfully tested (factor of two speedup with an error inferior to $1\%$) on heat transfer \cite{Bai2020} and fluid flow \cite{Bai2022} problems.
A similar approach relying on a more rigorous \textit{a posteriori} error estimator to switch between the HDM and ROM is introduced in \cite{Feng2021}.
A different AROM method developed in \cite{2020_Peherstorfer_AADEIM} uses the adaptive discrete empirical interpolation method (ADEIM) \cite{2015_Peherstorfer_ADEIM} and rank-one updates to adapt the reduced basis.
A comparison of AROMs relying on this approach and traditional ROMs can be found in \cite{2022_uy_aadeim}. In particular, the numerical experiments show that AROMs can be used in a predictive setting to model chemically reacting flow problems, whereas traditional ROMs completely fail to generate meaningful predictions.
The speedup factors achieved by these methods may seem at first very modest in comparison to the two-step ROM approach \cite{lui_wolf_2019}, but the offline phase cost is rarely discussed and too frequently only the online phase cost is taken into consideration. This is the case because online computational savings are assumed to be worth the offline cost. Since the amount of HDM data needed to produce accurate two-step ROMs strongly depends on the problem nonlinearity and parameter space size, highly nonlinear phenomena such as shock waves and turbulence may result in infeasible ROMs due to prohibitively costly offline phases.

In this work, we propose a training-free approach that combines local HDM and ROM solutions to cut down on costly full HDM solves. A dynamic relative reconstruction error strategy is developed to identify regions of the domain where the ROM is inaccurate and we locally solve the HDM in these regions. For problems containing spatial derivatives, states on neighboring cells are required to locally evolve the state using the HDM. We rely on the ROM solution when a neighboring cell is outside the sampled region.
Our approach allows the sampled region to adapt over time to avoid unnecessary HDM evaluations and improve robustness. 
Furthermore, our method relies on explicit spatial filtering combined with a residual-based error indicator to eliminate spurious oscillations that may appear after combining the solutions originating from different methods (e.g., some regions of the domain evolved using the HDM and others using the ROM). 
We refer to the solutions generated by this approach as hybrid snapshots.
These keys ingredients are novel contributions of this work, and important for stable and accurate prediction of shock-dominated flows as demonstrated on two canonical time-dependent compressible flow problems.

The remainder of this paper is organized as follows. 
In Section \ref{sec:arom}, we begin by introducing a general governing system of conservation laws and the high-dimensional modeling framework used to discretize it. 
Next, we introduce our hybrid snapshot approach, which involves: 1) the reduced basis approximation and partial HDM solutions, 2) a sampling procedure based on relative reconstruction error, and 3) low-pass spatial filters required to robustly mix solutions produced by different numerical methods.
We finish this section with a complete description of the algorithm and a discussion of important aspects of the method such as computational efficiency.
Section \ref{sec:num_exp} applies our adaptive framework to two compressible inviscid flow problems. The first is a compressible one-dimensional problem and is used to conduct a parametric study of the proposed method. The second is a considerably more complex two-dimensional problem.
Finally, Section \ref{sec:conclusions} highlights the main conclusions and discuss future research directions.

\section{Adaptive reduced-order models}
\label{sec:arom}

In this section, we introduce the general system of conservation laws that we aim to accelerate using our adaptive reduced-order model.
We begin by introducing the system of conservation laws (Section \ref{sec:claw}) and formulate a high-dimensional discretization (Section \ref{sec:hdm}).
Afterwards, we introduce our cost effective hybrid snapshot approach (Section \ref{sec:hybrid_snap}), which consists of the reduced basis approximation (Section \ref{sec:reduced_basis}), partial HDM solves (Section \ref{sec:partial_hdm}), relative reconstruction error (Section \ref{sec:rre}), and spatial low-pass filters (Section \ref{sec:filter}).

\subsection{System of conservation laws}
\label{sec:claw}

A general system of $c$ conservation laws, defined in a spatial domain $\Omega \subset \mathbb{R}^d$ over the time interval $\mathcal{T} = (0, T]$, takes the form 
\begin{equation}
    Q_{,t} + \nabla \cdot f (Q, \nabla Q) = h (Q, \nabla Q)
    , \; \; \; \; \; 
    Q(\cdot,0) = \mathring{Q} (\cdot),
    \label{eq:cons}
\end{equation}
where $f : \mathbb{R}^c \times \mathbb{R}^{c \times d} \rightarrow \mathbb{R}^{c \times d}$ is the flux function, $h : \mathbb{R}^c \times \mathbb{R}^{c \times d} \rightarrow \mathbb{R}^c$ is the source term, $\mathring{Q} : \Omega \rightarrow \mathbb{R}^c$ is the initial condition, and $Q(x,t)$ is the vector of conservative variables implicitly defined as the solution of Eq. \eqref{eq:cons} at $(x,t) \in \Omega \times \mathcal{T}$.

\subsection{High-dimensional model}
\label{sec:hdm}

The previous system of partial differential equations (PDEs) is discretized using a method of lines approach.
After spatial discretization, we have the following system of ordinary differential equations (ODEs)
\begin{equation}
    \frac{dq}{dt} = f^{*} (q,t)
%    \frac{dq}{dt} = F(q,t)
    \mbox{ ,}
    \label{eq:semi_disc}
\end{equation}
where $q(t) \in \mathbb{R}^N$ is our semi-discrete approximation to $Q (\cdot,t)$ implicitly defined as the solution of Eq. (\ref{eq:semi_disc}), $N$ is the number of degrees of freedom of the spatial discretization, and $f^{*}$ is the nonlinear function defining the spatial discretization of the inviscid and viscous fluxes. 

A time discretization method is required to solve Eq. \eqref{eq:semi_disc} numerically. In this work, the backward differentiation formulas (BDFs) are used. The $s$-order BDF scheme is written as
\begin{equation}
    \sum_{j=0}^{s} a_j q_{n+j} = \Delta t \beta f^{*} (q_{n+s},t_{n+s})
    \mbox{ ,}
    \label{eq:bdf}
\end{equation}
where $q_n \approx q (t_n)$, $\Delta t$ denotes the time step size, $t_n = t_1 + n \Delta t$, and coefficients $a_k$ and $\beta$ are such that the method is order $s$ and are normalized such that $a_s = 1$.
As can be noted from Eq. \eqref{eq:bdf}, BDF schemes are implicit and, thus, may require the solution of a nonlinear system of equations.

The fully discrete HDM is characterized by the following system of algebraic equations to be solved at each time instance $k \in [1, \ldots , N_t] $,
\begin{equation}
	R_k (q_k) \coloneqq q_k - F_k (q_k) = 0
	\mbox{ ,}
	\label{eq:hdm}
\end{equation}
where $R_k : \mathbb{R}^N \rightarrow \mathbb{R}^N$ is the nonlinear residual function and $F_k$ is a nonlinear function defined as
\begin{equation}
    F_k (q_k) \coloneqq \Delta t \beta f^{*} (q_k,t_k) - \sum_{j=0}^{s-1} a_j q_{k-s+j}
    \mbox{ .}
\end{equation}

\subsection{Hybrid snapshot approach}
\label{sec:hybrid_snap}

\begin{figure}[hbt!] % [H]
    \centering
    \includegraphics[width=.40\textwidth,trim={0mm 0mm 0mm 0mm},clip]{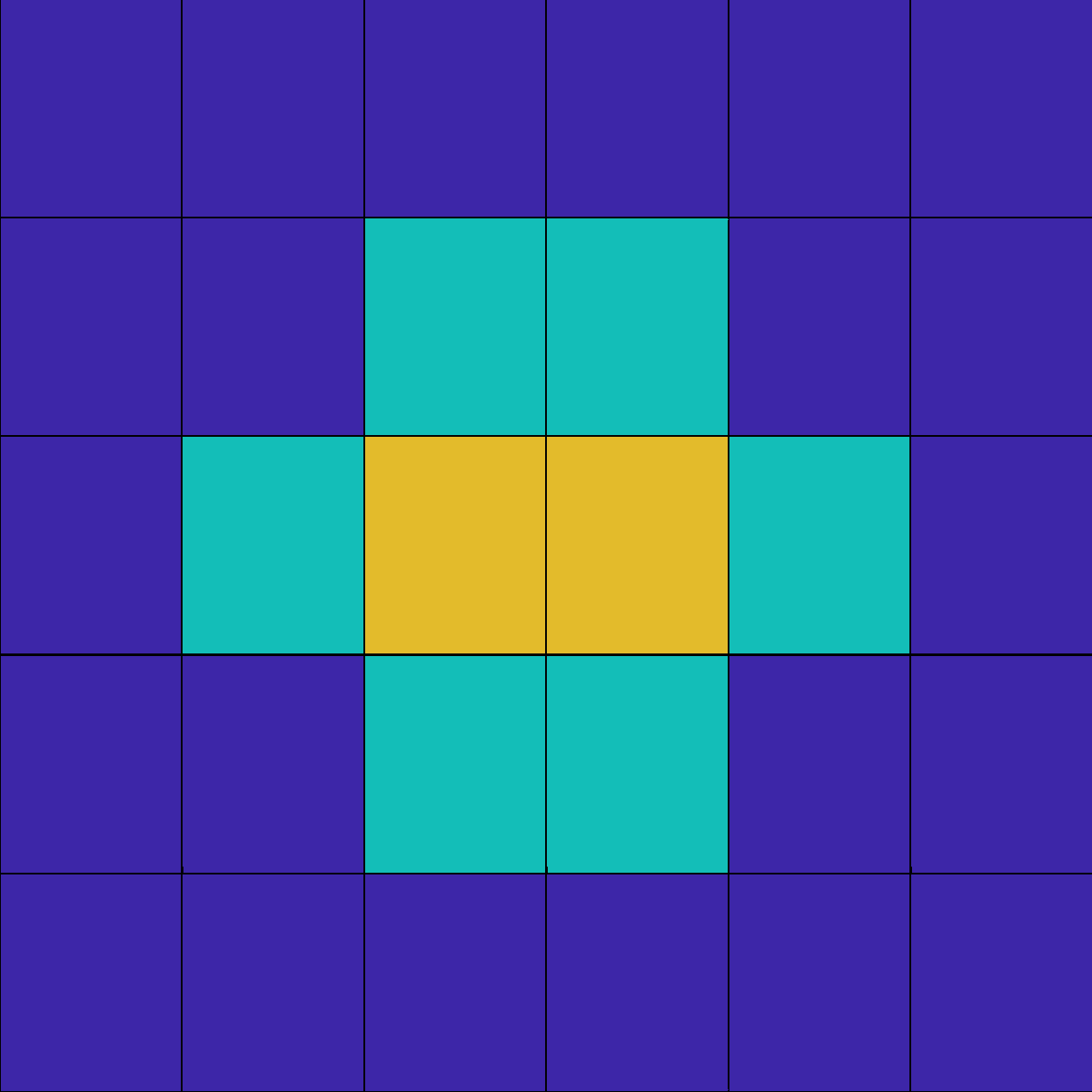}
    \caption{An example of mesh sampling corresponding to a first-order finite volume scheme.
    Cells sampled by $\hat{S}_k$ and $\tilde{S}_k$ are highlighted in yellow and teal respectively. Moreover, $\breve{S}_k$ samples both the blue and teal cells.}
    \label{fig:checkerboard}
\end{figure}

We are interested in obtaining an approximation $v_k \approx q_k$ that efficiently leverages local HDM information.
For this, consider the sampling points $\hat{s}_1^{(k)}, \ldots, \hat{s}_{n_s}^{(k)} \in \{1, \ldots, N \}$ and the corresponding sampling points matrix $\hat{S}_k = [e_{\hat{s}_1^{(k)}}, \dots , e_{\hat{s}_{n_s}^{(k)}}] \in \mathbb{R}^{N \times n_s}$. Here, $n_s$ is the number of indices retained from the original vector of size $N$ and $e_i$ denotes the vector with a $1$ in the $i$-th coordinate and $0$ elsewhere.
Let $\breve{S}_k \in \mathbb{R}^{N \times (N-n_s)}$ be the complementary sampling points matrix derived from points $\{ 1, \ldots , N \} \setminus \{ \hat{s}_1^{(k)}, \ldots , \hat{s}_{n_s}^{(k)} \}$ that have not been selected as sampling points.
We additionally consider sampling matrix $\tilde{S}_k \in \mathbb{R}^{N \times l}$ generated from the neighboring points $\{ \tilde{s}_1^{(k)}, \ldots , \tilde{s}_l^{(k)} \}$ needed to calculate the HDM flux function that are not already in $\{ \hat{s}_1^{(k)}, \ldots , \hat{s}_{n_s}^{(k)} \}$.
The sampling matrices are illustrated in Fig. \ref{fig:checkerboard} for the case of a first-order finite volume discretization.

With these definitions in place, we propose an approximation $v_k$ to the fully discrete HDM state $q_k$ where $v_k$ restricted to
the points in $\breve{S}_k$ use a traditional affine subspace approximation and $v_k$ restricted to the points in $\hat{S}_k$ are
defined as the solution of the HDM residual restricted to the $\hat{S}_k$ indices. That is, $v_k$ is defined such that
\begin{subequations}
    \begin{align}
        \breve{S}_k^\top v_k &= 
        \breve{S}_k^\top (\psi_k + \Phi_k y_k )
        \mbox{ ,}
        \label{eq:rom_approx}
        \\
        \hat{S}_k^\top v_k &= \hat{F}_k (\hat{S}_k^\top v_k, \tilde{S}_k^\top v_k)
        \mbox{ ,}
        \label{eq:partial_hdm}
    \end{align}
    \label{eq:hybrid}
\end{subequations}
where $\psi_k \in \mathbb{R}^N$ is a reference state, $\Phi_k \in \mathbb{R}^{N \times m}$ is a basis for a reduced subspace used to approximate the state $q_k$ at the sampling points $\breve{S}_k$, $y_k \in \mathbb{R}^m$ contains the corresponding reduced coordinates,   and $m$ denotes the dimension of the reduced subspace with $m \ll N$.
%\red{The function $\hat{F}_k : \mathbb{R}^{n_s} \times \mathbb{R}^l \rightarrow \mathbb{R}^{n_s}$ is the nonlinear partial residual defined as the restriction of the HDM residual $F_k$ to the indices sampled by $\hat{S}_k$.}
%
The function $\hat{F}_k : \mathbb{R}^{n_s} \times \mathbb{R}^l \rightarrow \mathbb{R}^{n_s}$ is defined as the restriction of the HDM nonlinear function $F_k$ to the indices sampled by $\hat{S}_k$.
Due to locality of the HDM discretization
scheme, the partial residual does not depend on the entire state; rather, it only depends on the restriction of the state to the
indices sampled by $\hat{S}_k$ and $\tilde{S}_k$. Mathematically, we write this as
\begin{equation} \label{eq:partial_hdm2}
    \hat{F}_k (\hat{v}, \tilde{v}) \coloneqq \hat{S}_k^\top F_k(\hat{S}_k \hat{v} + \tilde{S}_k \tilde{v})
    \mbox{ .}
\end{equation}
Evaluating $\hat{F}_k$ is cost effective provided $n_s \ll N$ because a relatively small number of entries of the HDM function are required.

\subsubsection{Reduced basis approximation}
\label{sec:reduced_basis}

We apply gappy POD \cite{Sirovich_gappy,Willcox_02_gappy} to compute the approximate HDM solution at the points corresponding to $\breve{S}_k$ (Eq. \ref{eq:rom_approx}).
Given a sampling matrix $P_k \in \mathbb{R}^{N \times n_p}$ constructed from points $\{ p_1^{(k)}, \ldots, p_{n_p}^{(k)} \} \subset \{\hat{s}_1^{(k)}, \ldots, \hat{s}_{n_s}^{(k)} \}$, the reduced coordinates $y_k$ are calculated as 
\begin{equation}
    y_k = (P_k^\top \Phi_k)^\dagger P_k^\top (v_k^{(J)} - \psi_k)
    \mbox{ ,}
    \label{eq:red_coord}
\end{equation}
where $v_k^{(J)}$ comes from the partial HDM solve, which is defined in Section \ref{sec:partial_hdm}.
The reduced basis $\Phi_k$ is constructed by compressing the deviations of the last $w$ snapshots from the reference state $\psi_{k-1}$, i.e., 
\begin{equation}
    % include q_bar?
    \Phi_k = {\tt POD}_m
    \left( \left[
    \gamma_{k-w} - \psi_{k-1}, \gamma_{k-w+1} - \psi_{k-1}, \ldots, \gamma_{k-2} - \psi_{k-1}, \gamma_{k-1} - \psi_{k-1}
    \right] \right)
    \mbox{ ,}
    \label{eq:pod}
\end{equation}
where $\gamma_k$ is either a HDM solution or hybrid snapshot (details deferred to Section \ref{sec:algo}) and ${\tt POD}_m : \mathbb{R}^{N \times w} \rightarrow \mathbb{R}^{N \times m}$ applies the thin SVD to the argument (snapshot matrix of size $N \times w$) and extracts the $m$ left singular vectors.
The sampling matrix $P_k$ is computed as
\begin{equation}
	% include q_bar?
	P_k = {\tt ODEIM}_{n_p} (\Phi_k)
	\mbox{ ,}
	\label{eq:odeim}
\end{equation}
where ${\tt ODEIM}_{n_p} : \mathbb{R}^{N \times m} \rightarrow \mathbb{R}^{N \times n_p}$ is the oversampling discrete empirical interpolation method (ODEIM) \cite{2020_Peherstorfer_ODEIM}, which is derived from the empirical interpolation method (EIM) \cite{maday_eim} and its discrete counterpart, the discrete empirical interpolation method (DEIM) \cite{Saifon01}.
As pointed out in \cite{2020_Peherstorfer_ODEIM}, oversampling ($m < n_p$) leads to more accurate linear-regression based approximations rather than interpolation ($m = n_p$).
Finally, the reference state is computed as
\begin{equation}
    \psi_k = \frac{1}{w} \sum_{j = k - w}^{k-1} \gamma_j
    \mbox{ .}
    \label{eq:mean}
\end{equation}
The reference state should be carefully chosen as it impacts accuracy and stability of the reduced bases approximation. 
In particular, our choice allows time-invariant Dirichlet boundary conditions to be automatically satisfied.

\begin{remark}
    Our SVD approach reconstructs the reduced basis from scratch every time the basis needs to be updated, which means all entries are updated. %In other words, all basis entries are updated.
    A different approach is adopted in \cite{2020_Peherstorfer_AADEIM}. In this case, the reduced-order basis is locally updated using the adaptive discrete empirical interpolation method (ADEIM) \cite{2015_Peherstorfer_ADEIM}.
    However, not providing any sort of correction outside the sampling points can lead to a potentially catastrophic loss of accuracy.
\end{remark}

\subsubsection{Partial high-dimensional model}
\label{sec:partial_hdm}

An estimate of $\tilde{v}_k \approx \tilde{S}_k^\top v_k$ is necessary in order to solve Eq. \ref{eq:partial_hdm} and, thus, obtain an approximate HDM solution at the points corresponding to $\hat{S}_k$. 
A straightforward choice is $\tilde{v}_k= \tilde{S}_k^\top\gamma_{k-1}$, i.e., lag the solution to the previous time step; however, this can lead to a lagged solution.
We attempt to obtain a more accurate evaluation of $\hat{S}_k^\top v_k$ through subiterations.
In this approach, solving the partial HDM solution at time step $k$ leads to the following iterations: for $j=1, \ldots , J$, solve
\begin{equation}
    \hat{v}_k^{(j)} = \hat{F}_k (\hat{v}_k^{(j)}, \tilde{v}_k^{(j)})
    % \mbox{ ,}
\end{equation}
for $\hat{v}_k^{(j)}$ and set
\begin{subequations}
    \begin{align}
        y_k^{(j)} &= (P_k^\top \Phi_k)^\dagger P_k^\top \hat{S}_k (\hat{v}_k^{(j)} - \psi_k)
        \mbox{ ,}
        \\
        \tilde{v}_k^{(j+1)} &= \tilde{S}_k^\top (\psi_k + \Phi_k y_k^{(j)})
    \end{align}
\end{subequations}
where $\tilde{v}_k^{(1)} = \tilde{S}_k^\top\gamma_{k-1}$ is the initial guess and $J$ is determined by the satisfaction of a convergence criterion.
Here, the algorithm is terminated when either 
\begin{equation}
    \| y_k^{(j+1)} - y_k^{(j)}\|_2 < \epsilon_y
\end{equation}
or $J = j_\mathrm{max}$, where $\epsilon_y \in \mathbb{R}_{>0}$ and $j_\mathrm{max} \in \mathbb{N}$ are user defined. 
In this work, we take $\epsilon_y = 10^{-4}$ and $j_\mathrm{max} = 10$ unless otherwise stated.

\begin{remark}
    For explicit time-marching methods, the right-hand size of Eq. \eqref{eq:hdm} can be directly computed because it only depends on the solution at previous time steps and, thus, no subiterations are necessary.
\end{remark}

\subsubsection{Relative reconstruction error}
\label{sec:rre}

The pointwise reconstruction error of approximating the state $\gamma_k$ in the reduced subspace is
\begin{equation}
    \varepsilon_{j}^{(k)} = (\gamma_k - \psi_k -  \Phi_k y_k)_j^2
    \mbox{ ,}
\end{equation}
where $y_k$ is given by Eq. \ref{eq:red_coord}.
Let $i_1, \ldots, i_N$ be an ordering such that
\begin{equation}
    \varepsilon_{i_1}^{(k)}  \geq \cdots  \geq \varepsilon_{i_N}^{(k)}
    \mbox{ .}
\end{equation}
% 
%At time step $k$, we pick the first $n_g$ indices $i_1 = g_1^{(k)},  \ldots, i_{n_g} = g_{n_g}^{(k)}$ as the sampling points to form $G_k$.
%
At time step $k$, we pick the first $n_g$ indices $i_1 = g_1^{(k)},  \ldots, i_{n_g} = g_{n_g}^{(k)}$ as the sampling points to form $G_k = [e_{g_1^{(k)}}, \dots , e_{g_{n_g}^{(k)}}] \in \mathbb{R}^{N \times n_g}$.
The number of sampling points $n_g$ is chosen according to the relative reconstruction error (RRE),
\begin{equation}
    \mbox{RRE} (n_g) = 
    \frac{\sum_{j=1}^{n_g} \varepsilon_{i_j}^{(k)}}{\sum_{j=1}^{N} \varepsilon_{i_j}^{(k)}}.
    \mbox{ ,}
    \label{eq:rre}
\end{equation}
In practice, we choose $n_g$ to be the smallest natural number such that $\mbox{RRE}(n_g) \geq \delta$.
This is done to identify the entries that concentrate most of the error.

Finally, the set of points forming sampling matrix $\hat{S}_k$ is defined as 
\begin{equation}
    \{ \hat{s}_1^{(k)}, \ldots , \hat{s}_{n_s}^{(k)} \} \coloneqq 
    \{ g_1^{(k)}, \ldots , g_{n_g}^{(k)} \} \cup \{ p_1^{(k)} , \ldots, p_{n_s}^{(k)} \}
    \mbox{ .}
\end{equation}
Sampling matrix $\hat{S}_k$ incorporates both the points that concentrate most of the reconstruction error and the ODEIM points needed to obtain reduced coefficients $y_k$.
Once we have $\hat{S}_k$, the other sampling matrices $\tilde{S}_k$ and $\breve{S}_k$ are straightforwardly obtained from the discrete stencil.

\begin{remark}
	Our approach is different from the sampling method presented in \cite{2020_Peherstorfer_AADEIM} in a few ways.
		First, we only take into consideration the last snapshot in the error evaluation. In \cite{2020_Peherstorfer_AADEIM}, all last $w$ snapshots are used.
		Second, the method introduced in \cite{2020_Peherstorfer_AADEIM} samples a fixed number of elements at all time instances. This can potentially lead to over- or under-sampling if the dynamically relevant region of the domain changes in size. On the other hand, our method fixes the error tolerance, which allows the number of sampled elements to adapt if required. 
		Third, sampling matrices are updated every $z$ time steps in \cite{2020_Peherstorfer_AADEIM}. In contrast, our method updates the sampling matrices more frequently. The increased accuracy provided by more frequent updates leads to smaller samples which, in turn, typically offsets the extra cost from the updates. This is particularly important for higher values of $z$.
		Finally, our method also incorporates the DEIM points. 
\end{remark}

\subsubsection{Spatial low-pass filters}
\label{sec:filter}

Spatial filtering is an operation commonly used to stabilize time-dependent fluid flow simulations \cite{Lele1992,VISBAL2002155,FALISSARD2013344} by eliminating high-wavenumber noise originating from, for example, mesh nonuniformity and nonlinear flow features.
Implicit filtering methods require the solution of a system of linear equations and have been used extensively in the solution of CFD problems \cite{Lele1992,VISBAL2002155}.
We avoid solving a system of linear equations by using the cheaper and easier to implement explicit filters. However, explicit filters require bigger stencils to obtain same order of accuracy which can be particularly problematic at boundaries.
Similar to standard CFD simulations, there is no guarantee that a hybrid solution $v_k$ combining entries from partial HDM and reduced basis solves is going to be smooth. To remove spurious oscillations that may develop, we apply one-dimensional explicit Shapiro filters \cite{FALISSARD2013344,shapiro_1970} to the solution.

For a hybrid solution, the nonlinear residual function equality defined by Eq. \refeq{eq:hdm} generally does not hold (i.e., $R_k (v_k) \neq 0$) and, thus, can be used as error indicator to selectively apply filters.
As a consequence, we can identify and restrict filtering to regions of the solution that lead to a non-negligible decrease of the residual function.
This avoids undesirable outcomes of filtering leading to an increase of the residual such as over-smoothing and filter-induced non-physical oscillations in the neighborhood of discontinuities and sharp gradients.
On the other hand, under-smoothing can also be detected by residual evaluation. In this case, the residual can be further reduced by successively applying filters of increasing order or repeated applications of the same filter operator \cite{shapiro_1970}.

We begin our local filtering procedure by computing the hybrid solution residual. Next, we apply the filter, recompute the residual, and discard the filtered solution on elements where the residual increases.
This procedure is repeated until the set of elements where the filtered solution is retained is empty, the absolute value of the element-wise residual decreases by less than $\epsilon_f \in \mathbb{R}_{>0}$, or the maximum number of iterations $j_\mathrm{max}^{f} \in \mathbb{N}$ is reached. In this work, we take $\epsilon_f = 10^{-2}$ and $j_\mathrm{max}^{f} = 10$ unless otherwise stated.

\begin{remark}
    As pointed out in \cite{VISBAL2002155}, multidimensional filtering can be performed by applying the one-dimensional filter in each coordinate direction.
\end{remark}

\begin{remark}
    Boundary condition treatment is usually not obvious and have been dealt with in different ways \cite{VISBAL2002155}.
    One approach is to use smaller, lower order stencils near the boundary, which decreases the global order of accuracy of the filter.
    Alternatively, decentered stencils maintaining the same order of accuracy as the centered stencil can be used. However, these need to be constructed in such a way that no frequency is amplified.
    In this work, for simplicity, the boundary values are obtained by using a zeroth-order extrapolation at the boundaries. 
\end{remark}

\begin{remark}
    Filtering is most commonly used on structured grids in combination with finite-difference methods. However, filtering can also be used on unstructured grids \cite{WOLFRAM201364}.
\end{remark}

\begin{remark}
	Residual evaluations and explicit filtering are operations that can be performed element-wise. Therefore, after the first iteration, this procedure becomes very cost effective if additional operations are only necessary at few elements.
\end{remark}

\begin{remark}
	Filtering stopping criteria tolerances were heuristically selected. Filters typically do not remove all frequency components above the given cutoff frequency. However, the benefits of successively applying the same filter rapidly diminish. 
	We observed that too tight residual tolerances or high number of iterations lead to a substantial increase in filtering cost (mostly from residual evaluations) with marginal accuracy increase. This can be particularly dramatic for higher-order filters.
\end{remark}

\subsection{General considerations, algorithm and computational efficiency}
\label{sec:algo}

\begin{figure}[hbt!]
    \centering
    \input{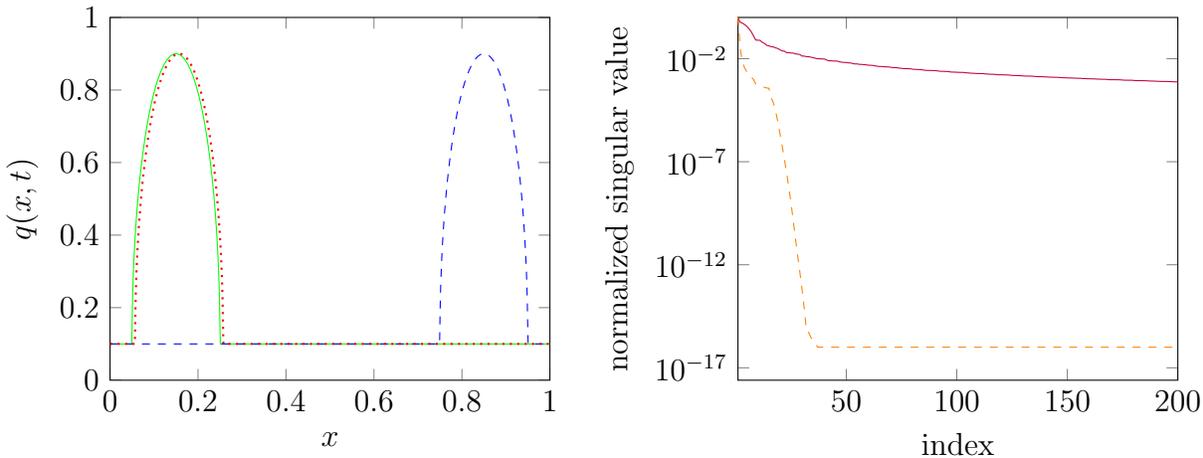}
    \caption{On the left, solutions of a scalar quantity at time $t = 0$ (\ref{line:fig2a_par1_t0}) and advected linearly at two different velocities $\mu_1$ (\ref{line:fig2a_par2_tf}) and  $\mu_2 = 100 \mu_1$ (\ref{line:fig2a_par1_tf}) for the same period of time.
    First and last snapshots of a scalar quantity advected linearly at two different velocities $(\mu_2 = 100 \mu_1)$ for the same period of time.
    On the right, normalized singular values for snapshots with global (\ref{line:fig2b_par1}) and local (\ref{line:fig2b_par2}) temporal structures.
    It can be noted that singular values of problems with local temporal structure decay orders of magnitudes faster compared to problems with global structure.
	}
    \label{fig:low_rank_struc}
\end{figure}

The proposed approach exploits the spatial locality of propagating coherent structures to derive efficient reduced-order models. As previously discussed, reduced-order modeling of convection-dominated problems is challenging because of the Kolmogorov barrier.
However, as pointed out in \cite{2020_Peherstorfer_AADEIM}, these problems have local low-rank structure: local trajectories have fast decaying singular values while the singular values of global trajectories decay slowly. 
The concept of local reduced bases for projection-based model reduction has also been exploited in other work \cite{Amsallem2012_local,Amsallem2015_local}.
A comparison of the trajectory of a scalar quantity advected linearly at two different velocities and their corresponding normalized singular values is illustrated in Fig. \ref{fig:low_rank_struc}.
As mentioned in Section \ref{sec:reduced_basis}, we construct the reduced basis by using the previous $w$ snapshots, where $w$ is chosen sufficiently small to ensure the subspace has a small dimension.

\begin{figure}[hbt!]
    \centering
    \subfloat[]{\includegraphics[width=.49\textwidth,trim={0mm 0mm 0mm 0mm},clip]{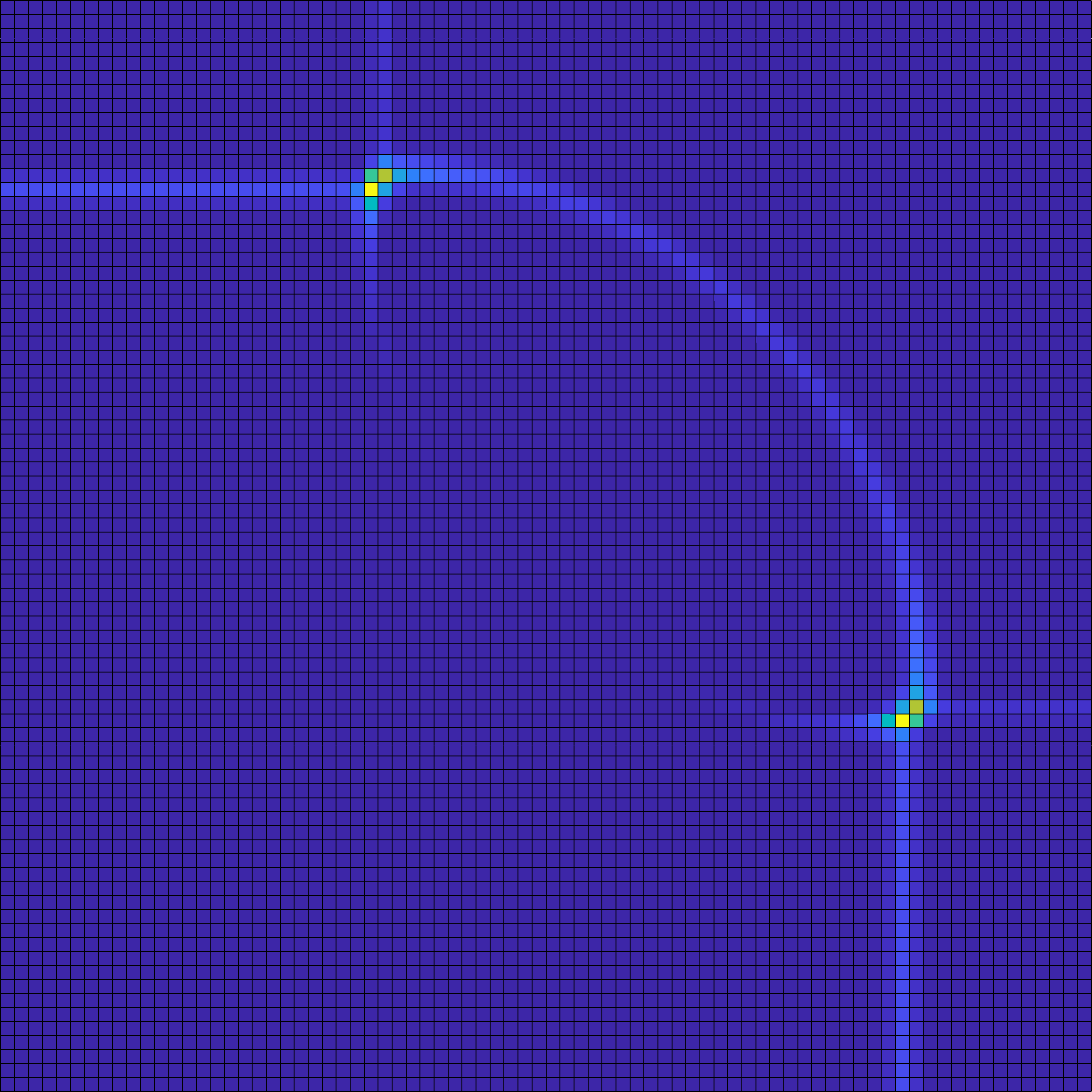}}
    ~
    \subfloat[]{\includegraphics[width=.49\textwidth,trim={0mm 0mm 0mm 0mm},clip]{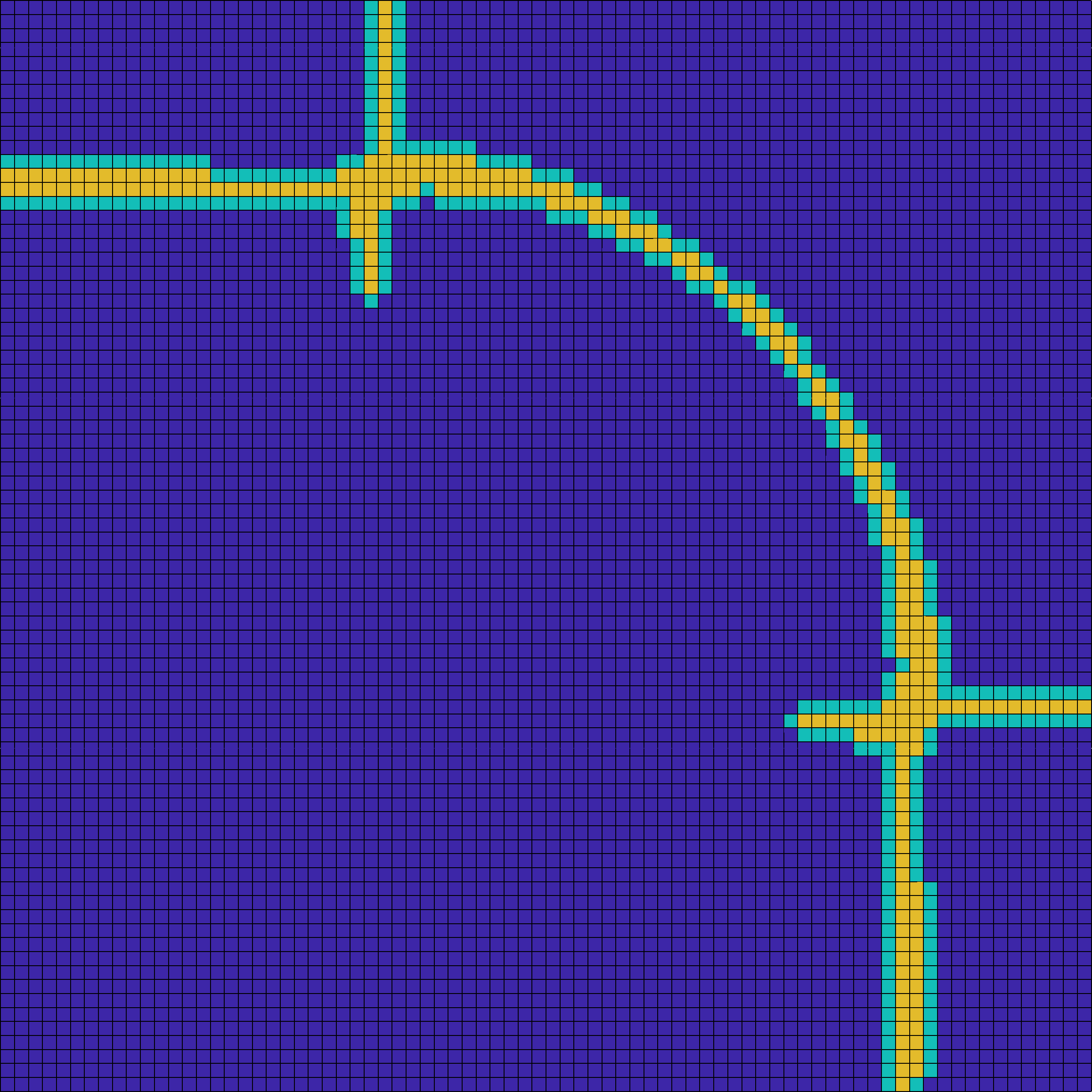}}
    \caption{Relative reconstruction error (left) and sampling points (right) of a problem with local spatial coherence.
    The entries corresponding to sampling matrices $\hat{S}_k$ and $\tilde{S}_k$ are highlighted in yellow and teal, respectively.
    The neighboring sampling matrix $\tilde{S}_k$ matches a first-order finite volume method, for example.
    }
    \label{fig:t17_res}
\end{figure}

Another important AROM ingredient is local spatial coherence. This feature leads to the RRE being concentrated at only a few components. In other words, the reduced basis is capable of providing an accurate approximation at most entries and, thus, more expensive HDM evaluations are only necessary at a small fraction of the components.
Figure \ref{fig:t17_res} illustrates an example of a problem where the RRE is concentrated in a few components only.
Entries where the RRE is small but nonzero will likely grow in time and result in an inaccurate solution. The proposed approach accounts for this by performing a full HDM solve every $z$ time steps.

\subsubsection{Algorithm}

Our AROM procedure is summarized in Algorithm \ref{alg:arom}.
The initial condition  is set at line \ref{alg_l:q0}.
The loop on line \ref{alg_l:for_timeStep} iterates over all time steps $k = 1, \ldots , N_t-1$.
The conditional statement on line \ref{alg_l:if_full_partial} chooses between a full (line \ref{alg_l:if_full}) or partial HDM solution (lines \ref{alg_l:if_partial_start}-\ref{alg_l:if_partial_end}). Initially, a full HDM solution is calculated for the first $w$ time steps. Afterwards, the second criterion ensures that a full HDM solution is going to take place every $z$ time steps.
A partial HDM computation takes place between lines \ref{alg_l:if_partial_start} and \ref{alg_l:if_partial_end}.
All other points are approximated via ODEIM (line \ref{alg_l:gappy}).
Line \ref{alg_l:filter} filters the hybrid snapshot originating from a partial HDM solution and RB reconstruction.
The conditional statement in line \ref{alg_l:if_rom} determines if the reduced basis and sampling points are computed. The first condition assures that the total number of snapshots is sufficient (i.e., at least $w$). The second condition checks if a full HDM evaluation is going to take place in the next time step. In this case, the reduced basis and sampling points are not necessary and, thus, do not need to be updated.
Finally, the reduced basis, sampling points and reference state are computed between lines \ref{alg_l:if_rom} and \ref{alg_l:if_rom_end}. The conditional on line \ref{alg_l:offset0} ensures the offset $\psi_k$ is available the first time the condition on line \ref{alg_l:if_rom} is satisfied.
The function ${\tt Neighbors}$ on line \ref{alg_l:neighbors} returns the sampling matrix $\tilde{S}_{k+1}$ generated from the neighboring points needed to calculate the HDM flux function that are not already sampled by $\hat{S}_{k+1}$ (Figure \ref{fig:checkerboard}). In addition, the set operations ${\tt union}$ and ${\tt setdiff}$ applied to sampling matrices are defined as the sampling matrix that results from the set operation applied to the index vector. That is, let $A,B \in \Rbb^{N\times N}$ be sampling matrices defined as $A = [e_{a_1},\dots,e_{a_s}]$ and
$B = [e_{b_1},\dots,e_{b_t}]$ from the index vectors $a\in\Nbb^s$, $b\in\Nbb^t$. Then,
\begin{equation}
C = {\tt union}(A,B), \qquad
D = {\tt setdiff}(A,B)
\end{equation}
are defined as the sampling matrices corresponding to the index vectors
$c = {\tt union}(a,b)$ and $d = {\tt setdiff}(a,b)$, respectively.

\begin{algorithm}[hbt!]
	\caption{Hybrid snapshot AROM}
	\begin{algorithmic}[1]
		\State Set $\gamma_0 = q_0$ \label{alg_l:q0}
		\For{$k = 1, \ldots, N_t$} \label{alg_l:for_timeStep}
		\If{$k+1 \leq w \OR {\tt mod}(k,z) = 0$} \label{alg_l:if_full_partial}
			\State Solve $\gamma_k = F_k (\gamma_k)$ for $\gamma_k$ \label{alg_l:if_full}
		\Else 
			\State $\tilde{v} = \tilde{S}_k^\top \gamma_{k-1}$ \label{alg_l:if_partial_start}
			\For{$j = 1, \ldots, J$}
			 	\State Solve $\hat{S}_k^\top \gamma_k = \hat{F}_k(\hat{S}_k^\top \gamma_k ,\tilde{v})$ for $\hat{S}_k^\top \gamma_k $
				\State $y_k = (P_k^\top \Phi_k)^\dagger P_k^\top (\gamma_k - \psi_k)$
				\State $\tilde{v} = \tilde{S}_k^\top (\psi_k + \Phi_k y_k)$
			\EndFor \label{alg_l:if_partial_end}
			\State $\breve{S}_k^\top \gamma_k = \breve{S}_k^\top (\psi_k + \Phi_k y_k)$ \label{alg_l:gappy}
			\State $\gamma_k = {\tt SpatialFilter} (\gamma_k)$ \label{alg_l:filter}
		\EndIf
		\If{$k = w-1$} \label{alg_l:offset0}
			\State $\displaystyle{\psi_{w-1} = \frac{1}{w}\sum_{j=0}^{w-1} \gamma_j}$
		\EndIf
		\If{$k \geq w-1 \AND {\tt mod}(k+1,z) \neq 0$} \label{alg_l:if_rom}
			\If{$k \geq w$}
				\State $G_{k+1}$ computed according to Section \ref{sec:rre}
			\Else
				\State $G_{k+1}$ = $\emptyset$
			\EndIf
			\State $\Phi_{k+1} = {\tt POD}_m (\left[
			\gamma_{k-w+1} - \psi_{k}, \gamma_{k-w+2} - \psi_{k}, \ldots, \gamma_{k-1} - \psi_{k}, \gamma_{k} - \psi_{k}
			\right])$
			\State $P_{k+1} = {\tt ODEIM}_{n_p}(\Phi_{k+1})$
			\State $\hat{S}_{k+1} = {\tt union} (G_{k+1}, P_{k+1})$
			\State $\breve{S}_{k+1} = {\tt setdiff} (I_N, \hat{S}_{k+1})$
			\State $\tilde{S}_{k+1} = {\tt Neighbors} ( \hat{S}_{k+1})$ \label{alg_l:neighbors}
			\State $\displaystyle{\psi_{k+1} = \frac{1}{w}\sum_{j=k-w+1}^{k} \gamma_j}$
		\EndIf \label{alg_l:if_rom_end}
		\EndFor
	\end{algorithmic}
	\label{alg:arom}
\end{algorithm}

\subsubsection{Computational efficiency}
\label{sec:complexity}

Our adaptive hybrid approach relies on $N$-dependent operations at every time step. 
A full HDM snapshot typically requires the solution of a nonlinear system by Newton's method, an iterative procedure that requires the solution of a linear system of equations at every time step. These large, sparse linear systems are usually solved with an iterative solver such as generalized minimal residual method (GMRES), which approximates the exact solve by a sequence of $\mathcal{O} (N^2)$ matrix-vector multiplications.
A hybrid snapshot computation (lines \ref{alg_l:if_partial_start}-\ref{alg_l:if_rom_end} of Algorithm \ref{alg:arom}) is going to require operations that at worst are log-linear.
For example, obtaining a reduced basis through a thin SVD and explicit filtering are algorithms that have linear complexity $\mathcal{O}(N)$.
A partial HDM iteration ($\mathcal{O} (n_s^2)$), selecting $n_p$ points with ODEIM $(\mathcal{O} (m^2 n_p^2))$ \cite{2020_Peherstorfer_ODEIM}, and computing the reduced coordinates $y_k$ through linear least squares $(\mathcal{O} (n_p m^2))$ are examples of operations independent of $N$. The RRE algorithm requires sorting the entries and, thus, is typically  $\mathcal{O} (N \log{N})$. While this sorting algorithm is the dominant term in terms of complexity, in practice it is not a bottleneck.

Let $t_H \in \mathbb{R}_{>0}$ and $t_R \in \mathbb{R}_{>0}$ be the average wall time required to compute a snapshot relying only on full HDM solutions and our adaptive approach, respectively.
Our AROM speedup $\mathcal{S}$ is defined in the following formula:
\begin{equation}
    \mathcal{S} \coloneqq \frac{t_H}{t_R}
    \mbox{ .}
    \label{eq:speedup}
\end{equation}
If the average sampling matrices are sufficiently small at all time steps such that the wall time required to compute a hybrid snapshot is negligible in comparison to a full HDM solution itis reasonable to assume $t_R \approx t_H / z$, which results in the following approximate speedup $\mathcal{S} \approx z$.
This implies the speedup of our approach is going to depend mainly on how often the full HDM must be solved if the sampling matrices remain reasonably small throughout the simulation.

\begin{remark}
    The complexity of obtaining a reduced basis through a thin SVD is $\mathcal{O} (N w^2)$. Therefore the number of snapshots used in the reconstruction $w$ is important to produce a small reduced basis but also a cost efficient construction.
    If necessary, reduced basis construction complexity can be reduced to $\mathcal{O} (N w^{\frac{1}{2}})$ by using fast SVD updates \cite{BRAND200620}.
\end{remark}

\begin{remark}
    In this work, we introduce a HDM that relies on BDF schemes for time-integration. 
    However, if an explicit scheme (e.g., Adams–Bashforth methods) was adopted instead, the computational complexity would be linear in $N$ as opposed to quadratic with an implicit scheme.
    For this class of ODE solvers, obtaining a cost efficient AROM can be considerably more challenging and problem dependent.
\end{remark}

\section{Numerical experiments}
\label{sec:num_exp}

In this section, we apply our adaptive method to solve two inviscid compressible flow problems.
We start by introducing the conservation laws, error functions and sampling average (Section \ref{sec:euler}).
The first test case is a canonical one-dimensional problem with known solution and is used to conduct a parametric study (Section \ref{sec:sod}) . For example, the impact of different filters and full HDM solve frequency are evaluated for this problem and serve as guideline for the next test case.
The second problem is two-dimensional and considerably more challenging (Section \ref{sec:implosion}).

\subsection{The Euler equations of gas dynamics}
\label{sec:euler}

We consider compressible inviscid flow through a domain $\Omega \subset \mathbb{R}^d$ with governing equations given by
\begin{subequations}
    \begin{equation}
        \frac{\partial \rho}{\partial t} + \frac{\partial}{\partial x_j} (\rho u_j) = 0
        \mbox{ }
    \end{equation}
    \begin{equation}
        \frac{\partial \rho u_i}{\partial t} + \frac{\partial}{\partial x_j} (\rho u_i u_j + P \delta_{ij}) = 0
        \mbox{ }
    \end{equation}
    \begin{equation}
        \frac{\partial (\rho E)}{\partial t} + \frac{\partial}{\partial x_j} ((\rho E + P) u_j) = 0
        \mbox{ ,}
    \end{equation}
    \label{eq:euler}
\end{subequations}
for $i = 1,\ldots, d$. The density of the fluid $\rho (\cdot, t) : \Omega \rightarrow \mathbb{R}_{>0}$, the fluid velocity $u(\cdot, t) \rightarrow \mathbb{R}^d$, and the total energy of the fluid $\rho E (\cdot , t) \rightarrow \mathbb{R}_{>0}$ are implicitly defined as the solution of \eqref{eq:euler}.
We assume the fluid follows the ideal gas law
\begin{equation}
    P = (\gamma - 1) \left(\rho E - \frac{\rho u_i u_i}{2}\right)
    \mbox{ ,}
\end{equation}
where $P (\cdot , t) \rightarrow \mathbb{R}_{>0}$ is the pressure of the fluid and $\gamma \in \mathbb{R}_{>0}$ is the ratio of specific heats.

We approximate the Euler equations using a finite volume method on a cartesian mesh. We employ a second-order monotonic upstream schemes for conservation laws (MUSCL) \cite{1979_vanLeer} approach with Roe flux \cite{roe_1981} and minmod limiter to spatially semi-discretize Eq. \eqref{eq:cons}. Afterwords, we integrate the resulting system of ODEs using a second-order BDF scheme defined by the coefficients $a_0 = 1/3$, $a_1 = - 4/3$, $a_2 = 1$ and $\beta = 2/3$.

In the following numerical experiments, the AROMs accuracy will be measured using the relative $L^1 (\Omega)$ error, defined as
\begin{equation}
    e_k \coloneqq 
    \frac{\int_{\Omega} \| \gamma_k (x) - q_k (x) \|_1 \,dV}{\int_{\Omega} \| q_k (x) \|_1 \,dV}
    \mbox{ .}
\end{equation}
% 
%\begin{equation}
%	e_k \coloneqq 
%	% 
%	\sqrt{\frac{\int_{\Omega} \| \gamma_k (x) - q_k (x) \|_2^2 \,dV}{\int_{\Omega} \| q_k (x) \|_2^2 \,dV}}
%	\mbox{ .}
%\end{equation}
% 
To access parametric performance, we also use the temporal mean of the relative error, defined as
\begin{equation}
    \bar{e} \coloneqq \frac{1}{N_t} \sum_{k=1}^{N_t} e_k
    \mbox{ .}
\end{equation}
Similarly, we define the average sampling as
\begin{equation}
    \bar{s} \coloneqq \frac{1}{N_t} \sum_{k=1}^{N_t} n_{\gamma _k}
    \mbox{ ,}
\end{equation}
where $n_{\gamma_k}$ is the number of entry points of snapshot $\gamma_k$ with its value directly computed by a HDM solve. For a snapshot $\gamma_k$ originating from partial and full HDM solves we have $n_{\gamma_k} = n_s$ and $n_{\gamma_k} = N$, respectively.
We define the average sampling of a hybrid snapshot as
\begin{equation}
	\bar{s}^{*} \coloneqq \frac{1}{| \Ical |} \sum_{k \in \Ical} n_{\gamma_k}
	\mbox{ ,}
\end{equation}
where $\Ical \subset \{ 1 , \ldots, N_t \}$ is the set of indices with a partial HDM solve. Lastly, we define the average ODEIM sampling as
\begin{equation}
	\bar{p} \coloneqq \frac{1}{| \Ical |} \sum_{k \in \Ical} (n_{p})_k
	\mbox{ .}
\end{equation}

\subsection{Sod's shock tube}
\label{sec:sod}

In this section we study our AROM method using the most canonical Riemann problem for the Euler equations, \textit{Sod's shock tube}.
We consider the one-dimensional ($d=1$) Euler equations in the domain $\Omega = (0,1)$ over the time interval $\mathcal{T} = (0, 0.2)$ with ratio of specific heats $\gamma = 1.4$ and initial condition, in terms of primitive variables, as
\begin{align}
	\rho (x,0) = 
	\begin{cases}
		1 & x \in [0, 0.5)\\
		0.125 & x \in [0.5, 1]
	\end{cases}
	,
	\quad
	u(x,0) = 0,
	\quad
	P(x,0) = 
	\begin{cases}
		1 & x \in [0, 0.5)\\
		0.1 & x \in [0.5, 1]
	\end{cases}
	\mbox{ .}
\end{align}
We use suitable boundary conditions from the initial condition. This is appropriate because the waves do not reach the boundary over the time interval of interest.

We partition the spatial domain into $N = 499$ cells of uniform width. We also equally partition the time domain into $N_t = 999$ time steps. 
The time step of implicit time marching methods is not limited by stability constraints that are typical of explicit methods. However, this does not imply that the time step can be arbitrarily large as it affects global accuracy. In this work, we performed a convergence study to chose a reasonable time step that is sufficiently small to ensure accurate solutions with steep discontinuity approximations.
The number of snapshots used in the reduced basis reconstruction is $w = 5$ and the number of POD modes used in the reconstruction is $m = 4$.
All hybrid solutions rely on the same reconstruction error threshold ($\delta = 0.80$). Moreover, we filter hybrid solutions by sequentially applying second-,fourth- and sixth-order filters. Lower-order filters are always applied first as they require a smaller number of filter passes to dissipate low-frequency noise.
These parameter values are used at all time steps unless otherwise stated.

\begin{figure}[hbt!]
	\centering
	\begin{tikzpicture}
\begin{groupplot} [
group style={group size = 2 by 2, horizontal sep = 1.8cm, vertical sep = 1.4cm}]
\nextgroupplot[width=.48\textwidth, xtick={2,4,6,8,10,12,14,16,18,20}, ytick={0,10,20,30,40,50,60}, xlabel={filter order}, ymax=60, xmax=10, ylabel={$\bar{s} (\%)$}, xmin=2, ymin=0]
\addplot [green, thick, mark options={solid, thin}, mark=o, mark size=3, mark repeat=0]
coordinates {
( 2.00000000e+00,  2.34717933e+00)
( 4.00000000e+00,  4.88251051e+00)
( 6.00000000e+00,  7.37045233e+00)
( 8.00000000e+00,  7.19274220e+00)
( 1.00000000e+01,  6.74776407e+00)
( 1.20000000e+01,  6.61282485e+00)
( 1.40000000e+01,  6.45961261e+00)
( 1.60000000e+01,  6.71061562e+00)
( 1.80000000e+01,  6.78973177e+00)
( 2.00000000e+01,  6.54896968e+00)};\label{line:fig_sod_filter_z1000_s}

\addplot [blue, mark options={solid, thin}, mark=diamond, mark size=3, mark repeat=0]
coordinates {
( 2.00000000e+00,  9.91452374e+00)
( 4.00000000e+00,  1.29746981e+01)
( 6.00000000e+00,  1.29223412e+01)
( 8.00000000e+00,  1.27101049e+01)
( 1.00000000e+01,  1.24392930e+01)
( 1.20000000e+01,  1.23476182e+01)
( 1.40000000e+01,  1.23472170e+01)
( 1.60000000e+01,  1.23022822e+01)
( 1.80000000e+01,  1.19590532e+01)
( 2.00000000e+01,  1.20653720e+01)};\label{line:fig_sod_filter_z15_s}

\addplot [red, thick, mark options={solid, thin}, mark=triangle, mark size=3, mark repeat=0]
coordinates {
( 2.00000000e+00,  2.42948823e+01)
( 4.00000000e+00,  2.48412657e+01)
( 6.00000000e+00,  2.46105437e+01)
( 8.00000000e+00,  2.44147614e+01)
( 1.00000000e+01,  2.43410669e+01)
( 1.20000000e+01,  2.41705857e+01)
( 1.40000000e+01,  2.41946819e+01)
( 1.60000000e+01,  2.40956864e+01)
( 1.80000000e+01,  2.41203851e+01)
( 2.00000000e+01,  2.40511082e+01)};\label{line:fig_sod_filter_z5_s}

\addplot [yellow, thick, mark options={solid, thin}, mark=square, mark size=3, mark repeat=0]
coordinates {
( 2.00000000e+00,  5.18731156e+01)
( 4.00000000e+00,  5.16735172e+01)
( 6.00000000e+00,  5.15523540e+01)
( 8.00000000e+00,  5.15304884e+01)
( 1.00000000e+01,  5.15034072e+01)
( 1.20000000e+01,  5.15445305e+01)
( 1.40000000e+01,  5.15431263e+01)
( 1.60000000e+01,  5.15690039e+01)
( 1.80000000e+01,  5.15509497e+01)
( 2.00000000e+01,  5.15136379e+01)};\label{line:fig_sod_filter_z2_s}

\nextgroupplot[width=.48\textwidth, xtick={2,4,6,8,10,12,14,16,18,20}, ytick={0,2,4,6,8,10}, xlabel={filter order}, ymax=8, xmax=10, ylabel={$\bar{s}^{*} (\%)$}, xmin=2, ymin=0]
\addplot [green, thick, mark options={solid, thin}, mark=o, mark size=3, mark repeat=0]
coordinates {
( 2.00000000e+00,  2.34717933e+00)
( 4.00000000e+00,  4.88251051e+00)
( 6.00000000e+00,  7.37045233e+00)
( 8.00000000e+00,  7.19274220e+00)
( 1.00000000e+01,  6.74776407e+00)
( 1.20000000e+01,  6.61282485e+00)
( 1.40000000e+01,  6.45961261e+00)
( 1.60000000e+01,  6.71061562e+00)
( 1.80000000e+01,  6.78973177e+00)
( 2.00000000e+01,  6.54896968e+00)};\label{line:fig_sod_filter_z1000_s_star}

\addplot [blue, thick, mark options={solid, thin}, mark=diamond, mark size=3, mark repeat=0]
coordinates {
( 2.00000000e+00,  3.24785708e+00)
( 4.00000000e+00,  6.30803148e+00)
( 6.00000000e+00,  6.25567451e+00)
( 8.00000000e+00,  6.04343823e+00)
( 1.00000000e+01,  5.77262633e+00)
( 1.20000000e+01,  5.68095149e+00)
( 1.40000000e+01,  5.68055029e+00)
( 1.60000000e+01,  5.63561558e+00)
( 1.80000000e+01,  5.29238657e+00)
( 2.00000000e+01,  5.39870532e+00)};\label{line:fig_sod_filter_z15_s_star}

\addplot [red, thick, mark options={solid, thin}, mark=triangle, mark size=3, mark repeat=0]
coordinates {
( 2.00000000e+00,  4.29488235e+00)
( 4.00000000e+00,  4.84126570e+00)
( 6.00000000e+00,  4.61054373e+00)
( 8.00000000e+00,  4.41476139e+00)
( 1.00000000e+01,  4.34106690e+00)
( 1.20000000e+01,  4.17058566e+00)
( 1.40000000e+01,  4.19468195e+00)
( 1.60000000e+01,  4.09568636e+00)
( 1.80000000e+01,  4.12038506e+00)
( 2.00000000e+01,  4.05110823e+00)};\label{line:fig_sod_filter_z5_s_star}

\addplot [yellow, thick, mark options={solid, thin}, mark=square, mark size=3, mark repeat=0]
coordinates {
( 2.00000000e+00,  1.87311560e+00)
( 4.00000000e+00,  1.67351720e+00)
( 6.00000000e+00,  1.55235396e+00)
( 8.00000000e+00,  1.53048840e+00)
( 1.00000000e+01,  1.50340721e+00)
( 1.20000000e+01,  1.54453050e+00)
( 1.40000000e+01,  1.54312629e+00)
( 1.60000000e+01,  1.56900387e+00)
( 1.80000000e+01,  1.55094975e+00)
( 2.00000000e+01,  1.51363789e+00)};\label{line:fig_sod_filter_z2_s_star}

\nextgroupplot[width=.48\textwidth, xtick={2,4,6,8,10,12,14,16,18,20}, ytick={0.01,1,100}, xlabel={filter order}, ymax=100, xmax=10, ylabel={$\bar{e} (\%)$}, xmin=2, ymode=log, ymin=0]
\addplot [green, thick, mark options={solid, thin}, mark=o, mark size=3, mark repeat=0]
coordinates {
( 2.00000000e+00,  1.79488387e+01)
( 4.00000000e+00,  2.95672681e+00)
( 6.00000000e+00,  3.15036096e-01)
( 8.00000000e+00,  3.01919527e-01)
( 1.00000000e+01,  2.35590967e-01)
( 1.20000000e+01,  2.08122161e-01)
( 1.40000000e+01,  1.83047051e-01)
( 1.60000000e+01,  1.78667163e-01)
( 1.80000000e+01,  1.48513317e-01)
( 2.00000000e+01,  1.58661989e-01)};\label{line:fig_sod_filter_z1000_e}

\addplot [blue, thick, mark options={solid, thin}, mark=diamond, mark size=3, mark repeat=0]
coordinates {
( 2.00000000e+00,  6.18083419e+00)
( 4.00000000e+00,  3.55662106e-01)
( 6.00000000e+00,  2.62879827e-01)
( 8.00000000e+00,  2.08469760e-01)
( 1.00000000e+01,  2.06419829e-01)
( 1.20000000e+01,  1.65241149e-01)
( 1.40000000e+01,  1.64761854e-01)
( 1.60000000e+01,  1.42273598e-01)
( 1.80000000e+01,  1.48756014e-01)
( 2.00000000e+01,  1.47393271e-01)};\label{line:fig_sod_filter_z15_e}

\addplot [red, thick, mark options={solid, thin}, mark=triangle, mark size=3, mark repeat=0]
coordinates {
( 2.00000000e+00,  5.07666468e-01)
( 4.00000000e+00,  1.51637694e-01)
( 6.00000000e+00,  1.21841155e-01)
( 8.00000000e+00,  1.05058580e-01)
( 1.00000000e+01,  1.08362016e-01)
( 1.20000000e+01,  9.25938930e-02)
( 1.40000000e+01,  9.24199230e-02)
( 1.60000000e+01,  7.83754964e-02)
( 1.80000000e+01,  8.88170737e-02)
( 2.00000000e+01,  8.24366935e-02)};\label{line:fig_sod_filter_z5_e}

\addplot [yellow, thick, mark options={solid, thin}, mark=square, mark size=3, mark repeat=0]
coordinates {
( 2.00000000e+00,  4.99908159e-02)
( 4.00000000e+00,  3.25336697e-02)
( 6.00000000e+00,  2.53204519e-02)
( 8.00000000e+00,  2.45241299e-02)
( 1.00000000e+01,  2.25930653e-02)
( 1.20000000e+01,  2.25004777e-02)
( 1.40000000e+01,  2.03074184e-02)
( 1.60000000e+01,  2.02340117e-02)
( 1.80000000e+01,  1.76479395e-02)
( 2.00000000e+01,  1.95558707e-02)};\label{line:fig_sod_filter_z2_e}

\end{groupplot}\end{tikzpicture}
	\caption{Time averages of relative sampling (top) and relative error (bottom) as a function of the filter order for full HDM frequency parameter $z=2$ (\ref{line:fig_sod_filter_z2_s}), $z=5$ (\ref{line:fig_sod_filter_z5_s}), $z=15$ (\ref{line:fig_sod_filter_z15_s}) and $z=\infty$ (\ref{line:fig_sod_filter_z1000_s}).
	A full HDM solve is equivalent to sampling all vector entries ($\bar{s} = 100\%$ and $\bar{s}^{*} = 0\%$).
	}
	\label{fig:t14_filter_error_sampling_}
\end{figure}
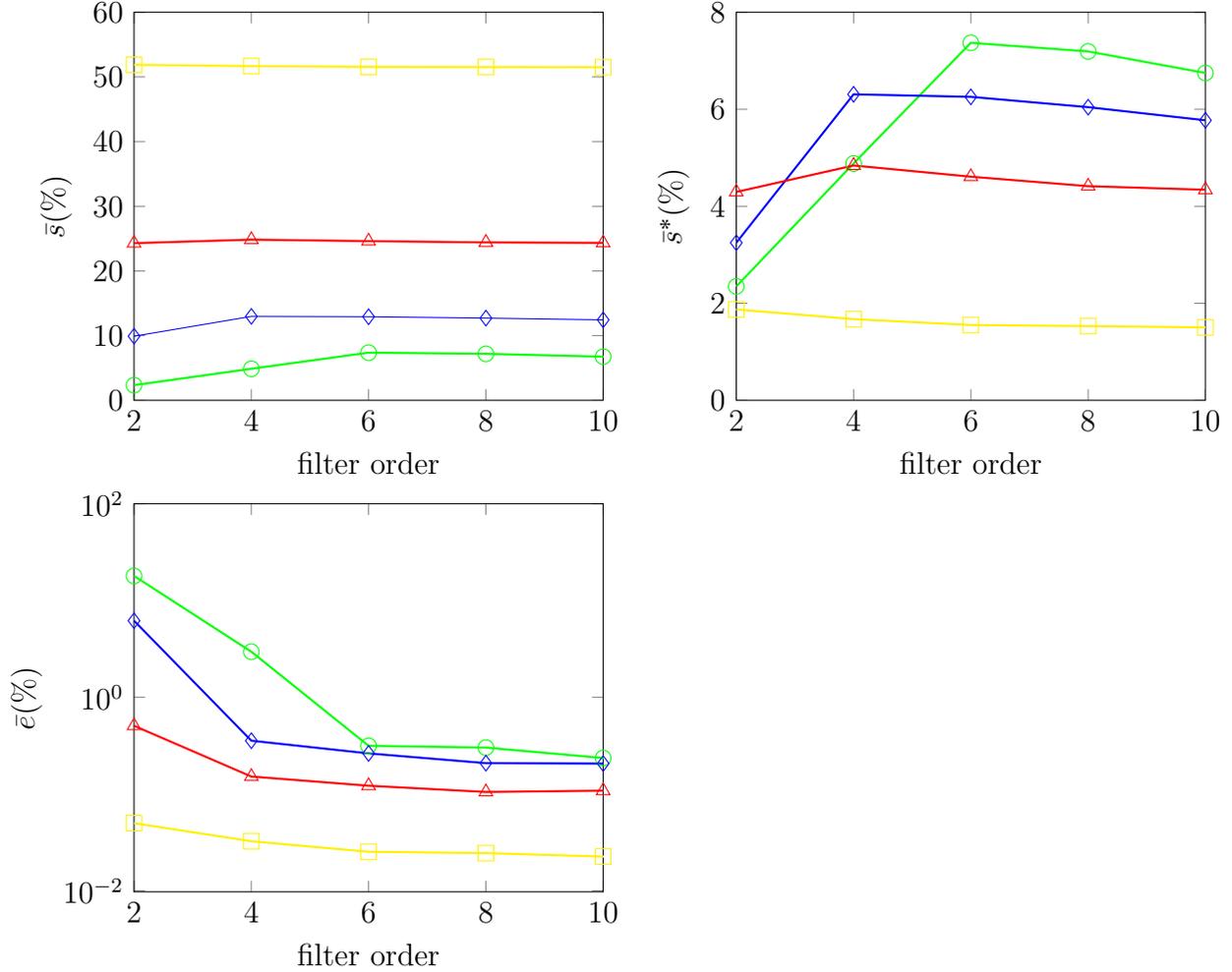

The effects of using filters of increasing accuracy on the hybrid solution can be observed in Fig. \ref{fig:t14_filter_error_sampling_}. 
High-order filters generally lead to more accurate solutions and particularly benefits simulations relying on low frequency full HDM solves the most.
For $z = 2$, all filters lead to accurate solutions with small sampling matrices. 
This is expected given that solving the full HDM every other time step results in most solution points being the result of HDM computations which in turn reduces the need for bigger sampling matrices at the partial HDM stage. Also, the HDM flux limiter inhibits the development of spurious oscillations. 
For all other values of $z$, only relying on lower-order filters leads to bigger errors and smaller sampling matrices. In fact, the under-damping of the hybrid solution, i.e., insufficient amount of viscosity to suppress all spurious oscillations, causes the RRE to be less equally distributed among the entries which in turn leads to smaller sampling matrices.

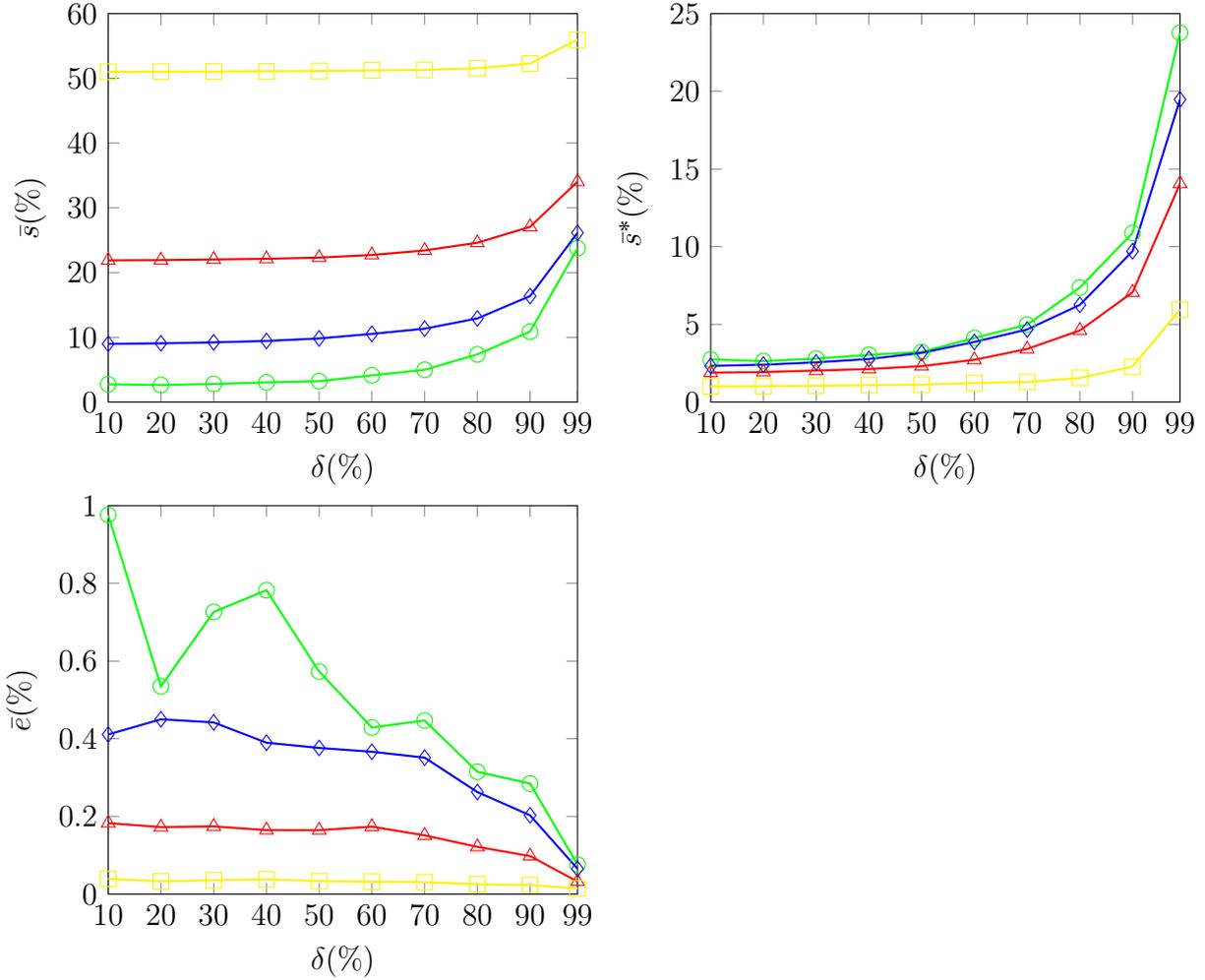
\begin{figure}[hbt!]
	\centering
	\begin{tikzpicture}
\begin{groupplot} [
group style={group size = 2 by 2, horizontal sep = 1.8cm, vertical sep = 1.4cm}]
\nextgroupplot[width=.48\textwidth, xtick={10,20,30,40,50,60,70,80,90,99}, ytick={0,10,20,30,40,50,60}, xlabel={$\delta (\%)$}, ymax=60, xmax=99, ylabel={$\bar{s} (\%)$}, xmin=10, ymin=0]
\addplot [green, thick, mark options={solid, thin}, mark=o, mark size=3, mark repeat=0]
coordinates {
( 1.00000000e+01,  2.74798093e+00)
( 2.00000000e+01,  2.63713800e+00)
( 3.00000000e+01,  2.80059116e+00)
( 4.00000000e+01,  3.03814041e+00)
( 5.00000000e+01,  3.22649307e+00)
( 6.00000000e+01,  4.12669025e+00)
( 7.00000000e+01,  4.98331332e+00)
( 8.00000000e+01,  7.37045233e+00)
( 9.00000000e+01,  1.08806792e+01)
( 9.90000000e+01,  2.37617520e+01)};\label{line:fig_sod_rre_z1000_s}

\addplot [blue, thick, mark options={solid, thin}, mark=diamond, mark size=3, mark repeat=0]
coordinates {
( 1.00000000e+01,  8.99637112e+00)
( 2.00000000e+01,  9.07259965e+00)
( 3.00000000e+01,  9.22485612e+00)
( 4.00000000e+01,  9.43869721e+00)
( 5.00000000e+01,  9.83067236e+00)
( 6.00000000e+01,  1.05363881e+01)
( 7.00000000e+01,  1.13385931e+01)
( 8.00000000e+01,  1.29223412e+01)
( 9.00000000e+01,  1.63744907e+01)
( 9.90000000e+01,  2.61373598e+01)};\label{line:fig_sod_rre_z15_s}

\addplot [red, thick, mark options={solid, thin}, mark=triangle, mark size=3, mark repeat=0]
coordinates {
( 1.00000000e+01,  2.18936872e+01)
( 2.00000000e+01,  2.19320404e+01)
( 3.00000000e+01,  2.20294296e+01)
( 4.00000000e+01,  2.21274212e+01)
( 5.00000000e+01,  2.23137658e+01)
( 6.00000000e+01,  2.27217963e+01)
( 7.00000000e+01,  2.34175766e+01)
( 8.00000000e+01,  2.46105437e+01)
( 9.00000000e+01,  2.70617387e+01)
( 9.90000000e+01,  3.40562889e+01)};\label{line:fig_sod_rre_z5_s}

\addplot [yellow, thick, mark options={solid, thin}, mark=square, mark size=3, mark repeat=0]
coordinates {
( 1.00000000e+01,  5.10025055e+01)
( 2.00000000e+01,  5.10227663e+01)
( 3.00000000e+01,  5.10470390e+01)
( 4.00000000e+01,  5.10899677e+01)
( 5.00000000e+01,  5.11178513e+01)
( 6.00000000e+01,  5.12219634e+01)
( 7.00000000e+01,  5.12971890e+01)
( 8.00000000e+01,  5.15523540e+01)
( 9.00000000e+01,  5.22879593e+01)
( 9.90000000e+01,  5.59543511e+01)};\label{line:fig_sod_rre_z2_s}

\nextgroupplot[width=.48\textwidth, xtick={10,20,30,40,50,60,70,80,90,99}, xlabel={$\delta (\%)$}, ymax=25, xmax=99, ylabel={$\bar{s}^{*} (\%)$}, xmin=10, ymin=0]
\addplot [green, thick, mark options={solid, thin}, mark=o, mark size=3, mark repeat=0]
coordinates {
( 1.00000000e+01,  2.74798093e+00)
( 2.00000000e+01,  2.63713800e+00)
( 3.00000000e+01,  2.80059116e+00)
( 4.00000000e+01,  3.03814041e+00)
( 5.00000000e+01,  3.22649307e+00)
( 6.00000000e+01,  4.12669025e+00)
( 7.00000000e+01,  4.98331332e+00)
( 8.00000000e+01,  7.37045233e+00)
( 9.00000000e+01,  1.08806792e+01)
( 9.90000000e+01,  2.37617520e+01)};\label{line:fig_sod_rre_z1000_s_star}

\addplot [blue, thick, mark options={solid, thin}, mark=diamond, mark size=3, mark repeat=0]
coordinates {
( 1.00000000e+01,  2.32970445e+00)
( 2.00000000e+01,  2.40593299e+00)
( 3.00000000e+01,  2.55818945e+00)
( 4.00000000e+01,  2.77203055e+00)
( 5.00000000e+01,  3.16400569e+00)
( 6.00000000e+01,  3.86972142e+00)
( 7.00000000e+01,  4.67192644e+00)
( 8.00000000e+01,  6.25567451e+00)
( 9.00000000e+01,  9.70782406e+00)
( 9.90000000e+01,  1.94706931e+01)};\label{line:fig_sod_rre_z15_s_star}

\addplot [red, thick, mark options={solid, thin}, mark=triangle, mark size=3, mark repeat=0]
coordinates {
( 1.00000000e+01,  1.89368717e+00)
( 2.00000000e+01,  1.93204043e+00)
( 3.00000000e+01,  2.02942960e+00)
( 4.00000000e+01,  2.12742118e+00)
( 5.00000000e+01,  2.31376581e+00)
( 6.00000000e+01,  2.72179630e+00)
( 7.00000000e+01,  3.41757664e+00)
( 8.00000000e+01,  4.61054373e+00)
( 9.00000000e+01,  7.06173871e+00)
( 9.90000000e+01,  1.40562889e+01)};\label{line:fig_sod_rre_z5_s_star}

\addplot [yellow, thick, mark options={solid, thin}, mark=square, mark size=3, mark repeat=0]
coordinates {
( 1.00000000e+01,  1.00250551e+00)
( 2.00000000e+01,  1.02276625e+00)
( 3.00000000e+01,  1.04703902e+00)
( 4.00000000e+01,  1.08996772e+00)
( 5.00000000e+01,  1.11785132e+00)
( 6.00000000e+01,  1.22196345e+00)
( 7.00000000e+01,  1.29718897e+00)
( 8.00000000e+01,  1.55235396e+00)
( 9.00000000e+01,  2.28795930e+00)
( 9.90000000e+01,  5.95435114e+00)};\label{line:fig_sod_rre_z2_s_star}

\nextgroupplot[width=.48\textwidth, xtick={10,20,30,40,50,60,70,80,90,99}, ytick={0,.2,.4,.6,.8,1}, xlabel={$\delta (\%)$}, ymax=1, xmax=99, ylabel={$\bar{e} (\%)$}, xmin=10, ymin=0]
\addplot [green, thick, mark options={solid, thin}, mark=o, mark size=3, mark repeat=0]
coordinates {
( 1.00000000e+01,  9.76055365e-01)
( 2.00000000e+01,  5.34950305e-01)
( 3.00000000e+01,  7.26400355e-01)
( 4.00000000e+01,  7.82271371e-01)
( 5.00000000e+01,  5.72991068e-01)
( 6.00000000e+01,  4.28847162e-01)
( 7.00000000e+01,  4.46972999e-01)
( 8.00000000e+01,  3.15036096e-01)
( 9.00000000e+01,  2.84744771e-01)
( 9.90000000e+01,  7.59182349e-02)};\label{line:fig_sod_rre_z1000_e}

\addplot [blue, thick, mark options={solid, thin}, mark=diamond, mark size=3, mark repeat=0]
coordinates {
( 1.00000000e+01,  4.10960956e-01)
( 2.00000000e+01,  4.50015622e-01)
( 3.00000000e+01,  4.42012366e-01)
( 4.00000000e+01,  3.89793428e-01)
( 5.00000000e+01,  3.76281295e-01)
( 6.00000000e+01,  3.66489955e-01)
( 7.00000000e+01,  3.51370251e-01)
( 8.00000000e+01,  2.62879827e-01)
( 9.00000000e+01,  2.02961719e-01)
( 9.90000000e+01,  6.53379712e-02)};\label{line:fig_sod_rre_z15_e}

\addplot [red, thick, mark options={solid, thin}, mark=triangle, mark size=3, mark repeat=0]
coordinates {
( 1.00000000e+01,  1.82923317e-01)
( 2.00000000e+01,  1.72623643e-01)
( 3.00000000e+01,  1.74368781e-01)
( 4.00000000e+01,  1.64922669e-01)
( 5.00000000e+01,  1.64727882e-01)
( 6.00000000e+01,  1.73828301e-01)
( 7.00000000e+01,  1.51315252e-01)
( 8.00000000e+01,  1.21841155e-01)
( 9.00000000e+01,  9.80864955e-02)
( 9.90000000e+01,  3.23462537e-02)};\label{line:fig_sod_rre_z5_e}

\addplot [yellow, thick, mark options={solid, thin}, mark=square, mark size=3, mark repeat=0]
coordinates {
( 1.00000000e+01,  3.92392276e-02)
( 2.00000000e+01,  3.29878884e-02)
( 3.00000000e+01,  3.57042713e-02)
( 4.00000000e+01,  3.80662731e-02)
( 5.00000000e+01,  3.37185569e-02)
( 6.00000000e+01,  3.20787920e-02)
( 7.00000000e+01,  3.10009759e-02)
( 8.00000000e+01,  2.53204519e-02)
( 9.00000000e+01,  2.33322182e-02)
( 9.90000000e+01,  1.47105679e-02)};\label{line:fig_sod_rre_z2_e}

\end{groupplot}\end{tikzpicture}
	\caption{Time averages of relative sampling (top) and relative error (bottom) as a function of relative reconstruction error tolerance $\delta$ for full HDM frequency parameter $z=2$ (\ref{line:fig_sod_rre_z2_s}), $z=5$ (\ref{line:fig_sod_rre_z5_s}), $z=15$ (\ref{line:fig_sod_rre_z15_s}) and $z=\infty$ (\ref{line:fig_sod_rre_z1000_s}).
	A full HDM solve is equivalent to sampling all vector entries ($\bar{s} = 100\%$ and $\bar{s}^{*} = 0\%$).
	}
	\label{fig:t14_rre_study_}
\end{figure}

The implication of different values of RRE tolerance $\delta$ can be observed in Fig. \ref{fig:t14_rre_study_}.
For $z = 2$, the error variation is negligible for the range of RRE tolerances considered. In regards to the time average sampling, it remains visually constant for most values of $\delta$ but abruptly increases for tighter tolerances. This shows, for this case in particular, that a smaller sampling matrix is enough to generate accurate AROMs.
For all other full HDM solve frequencies, accuracy can be more significantly improved by the use of tighter RRE tolerances. This is particularly substantial when full HDM solves are only performed on the initial training stage ($z = \infty$). On the other hand, accuracy comes at a price as bigger sampling matrices become necessary.
Moreover, increasing the RRE tolerance  did not lead to the time average error to monotonically decrease. One possible explanation is that adding just a few sampling points could add noise to solution. In general, having more solution points originating from a partial HDM solution leads to a more accurate AROM. However, this could introduce undesirable higher frequency structures, especially if the points are sparsely distributed.

\begin{figure}[hbt!]
	\centering
	\subfloat[$z=2$]{\includegraphics[width=.49\textwidth,trim={0mm 0mm 0mm 0mm},clip]{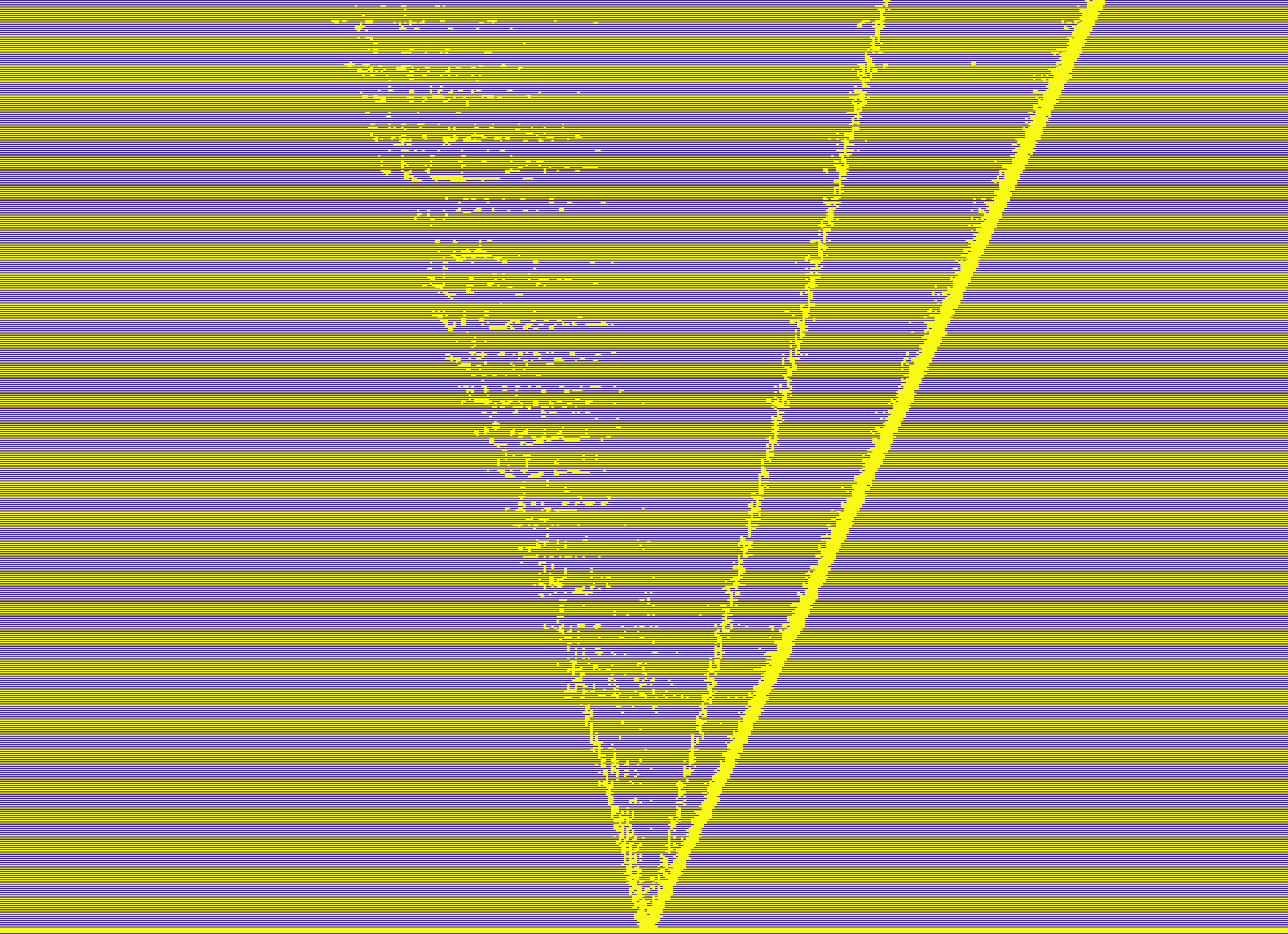}}
	~
	\subfloat[$z=5$]{\includegraphics[width=.49\textwidth,trim={0mm 0mm 0mm 0mm},clip]{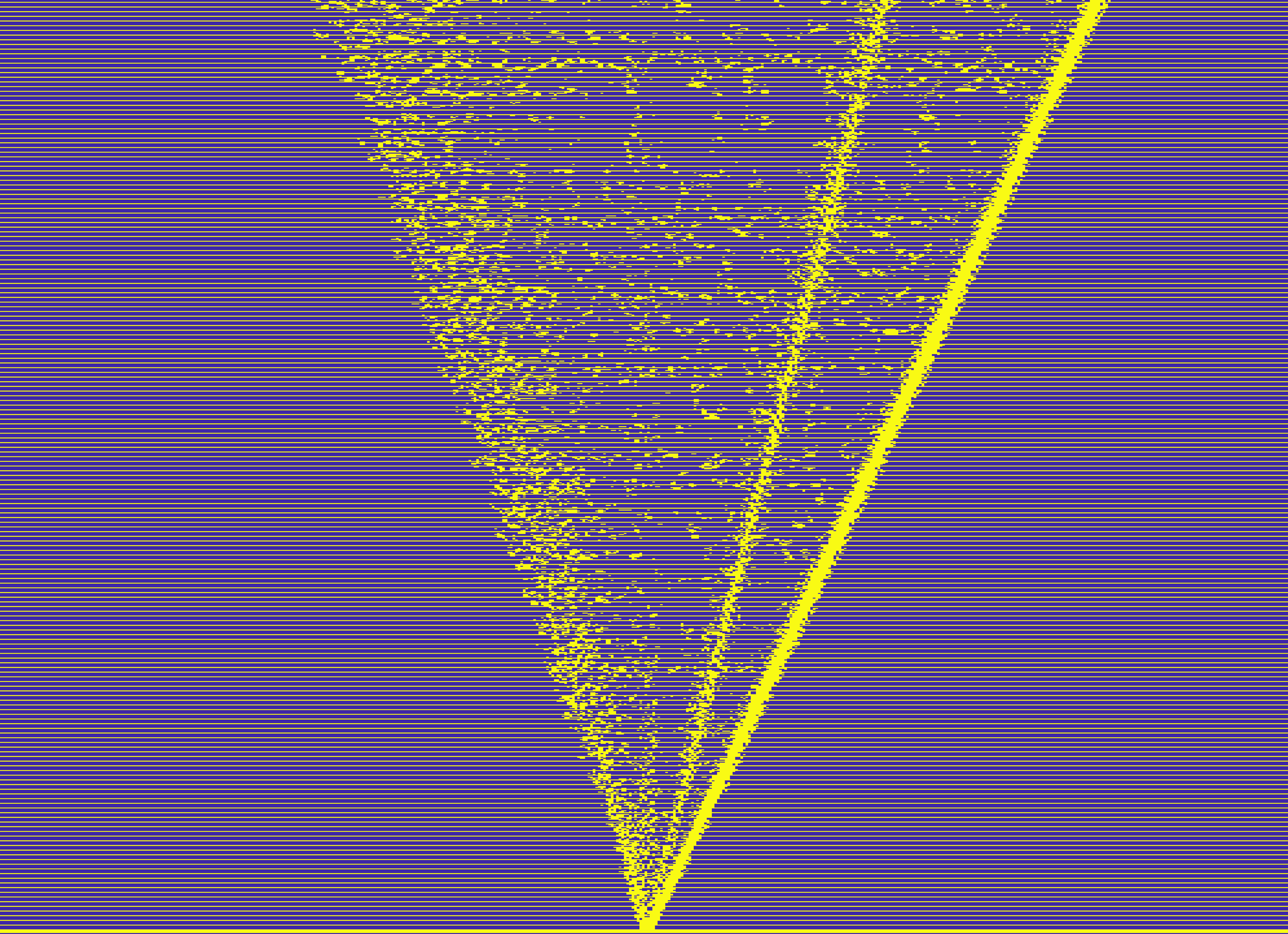}}
	\\
	\subfloat[$z=15$]{\includegraphics[width=.49\textwidth,trim={0mm 0mm 0mm 0mm},clip]{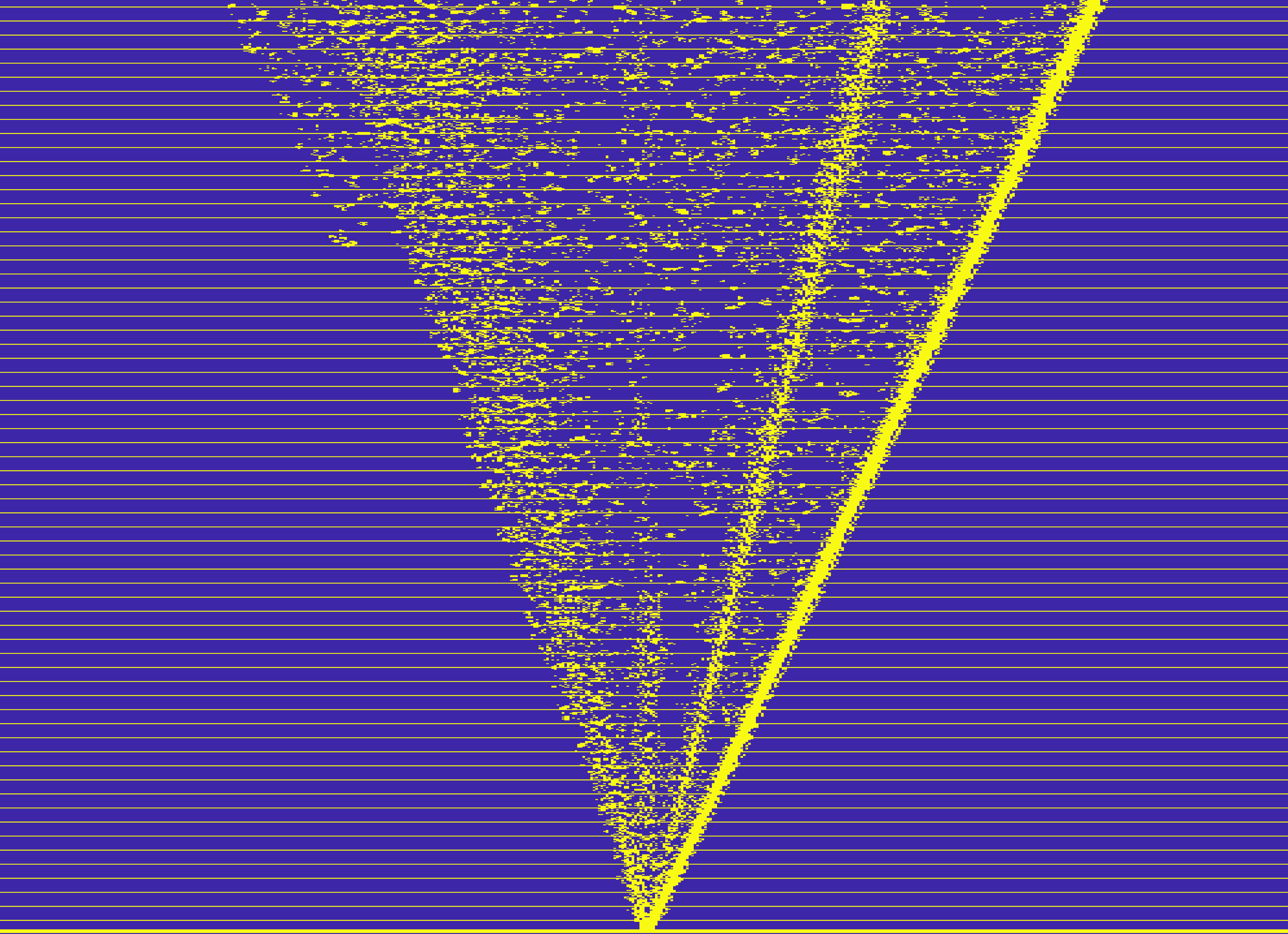}}
	~
	\subfloat[$z=\infty$]{\includegraphics[width=.49\textwidth,trim={0mm 0mm 0mm 0mm},clip]{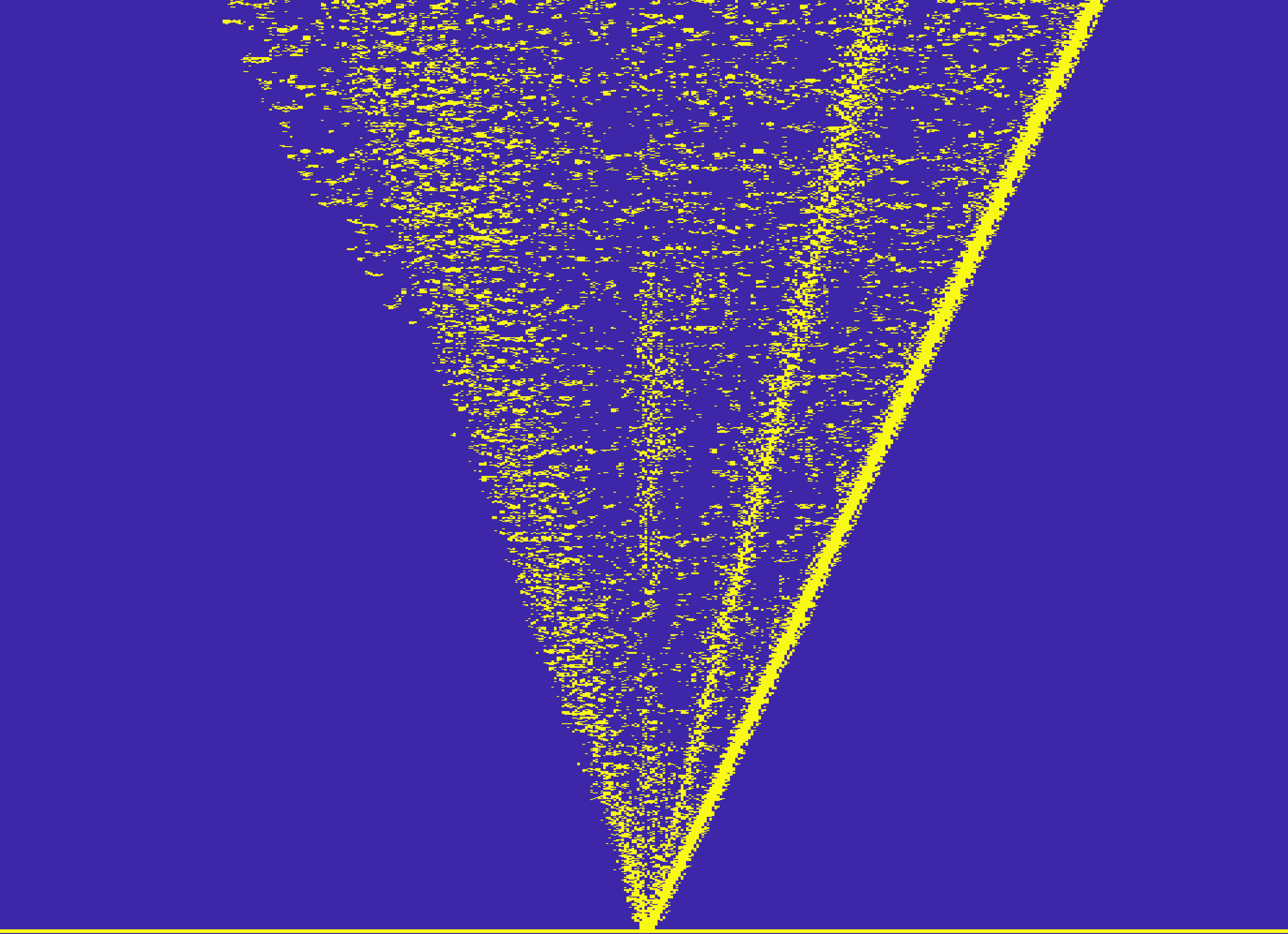}}
	\caption{Space-time snapshots of sampling points selected by matrix $\hat{S}$ (in yellow) for AROMs with different full HDM solve frequency $z$.}
	\label{fig:t14_map_}
\end{figure}

\begin{figure}[hbt!]
	\centering
	\input{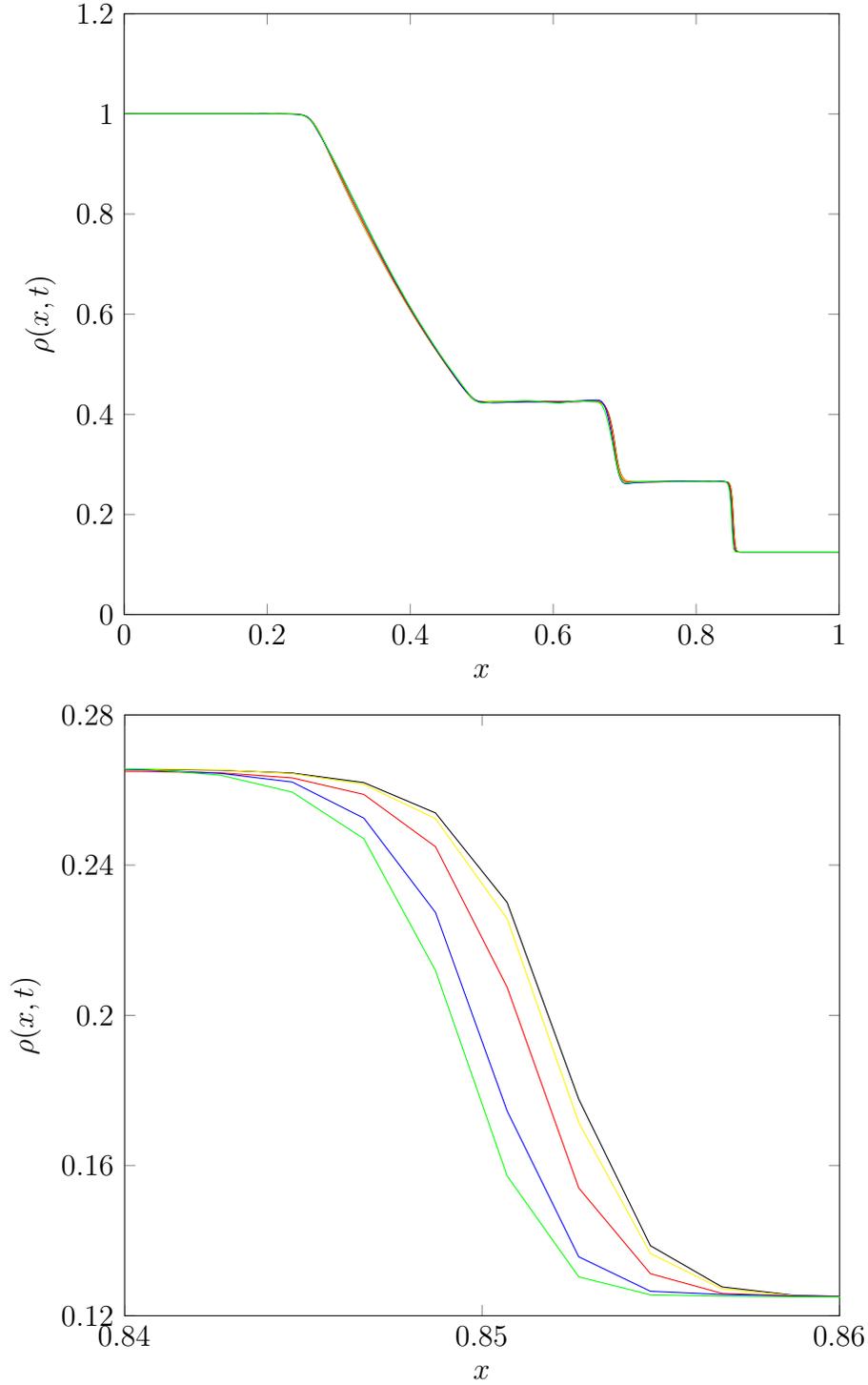}
	\caption{The HDM (\ref{line:fig_sod_q_hdm}) solution (density) at $k = N_t$ and the corresponding AROMs for $z=2$ (\ref{line:fig_sod_q_z2}), $z=5$ (\ref{line:fig_sod_q_z5}), $z=15$ (\ref{line:fig_sod_q_z15}) and $z = \infty$ (\ref{line:fig_sod_q_z1000}).}
	\label{fig:t14_snapshot_}
\end{figure}

\begin{figure}[hbt!]
	\centering
	\input{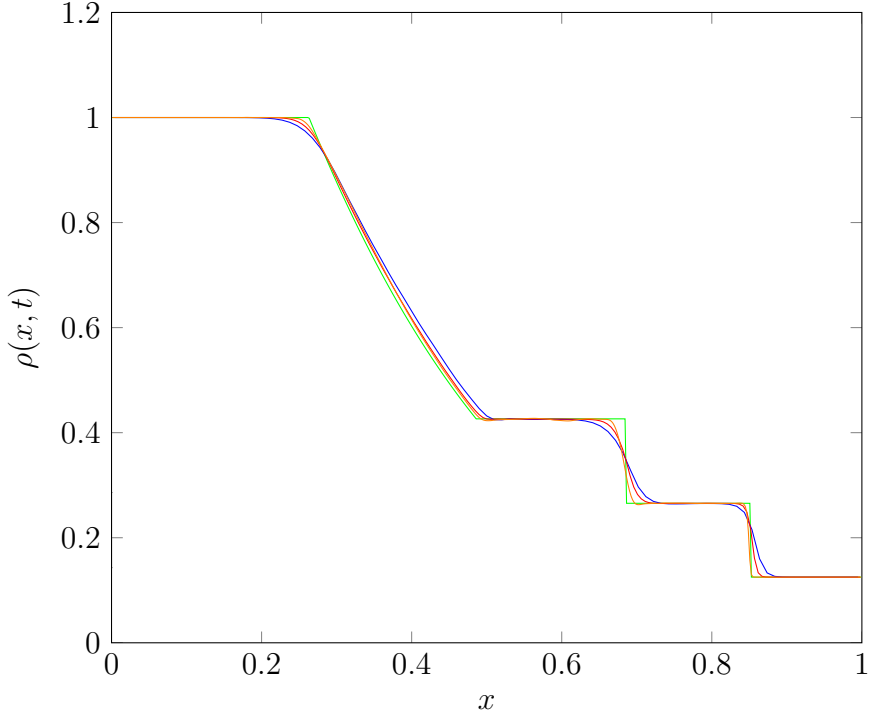}
	\caption{The exact (\ref{line:fig_sod_exact}) solution (density) at $k = N_t$, and the corresponding HDM for $N = 99$ (\ref{line:fig_sod_hdm100}), $N = 199$ (\ref{line:fig_sod_hdm200}), and the corresponding AROM for $z = \infty$ (\ref{line:fig_sod_arom500}).}
	\label{fig:t14_coarse}
\end{figure}

\begin{figure}[hbt!]
	\centering
	\subfloat[]{\includegraphics[width=.49\textwidth,trim={0mm 0mm 0mm 0mm},clip]{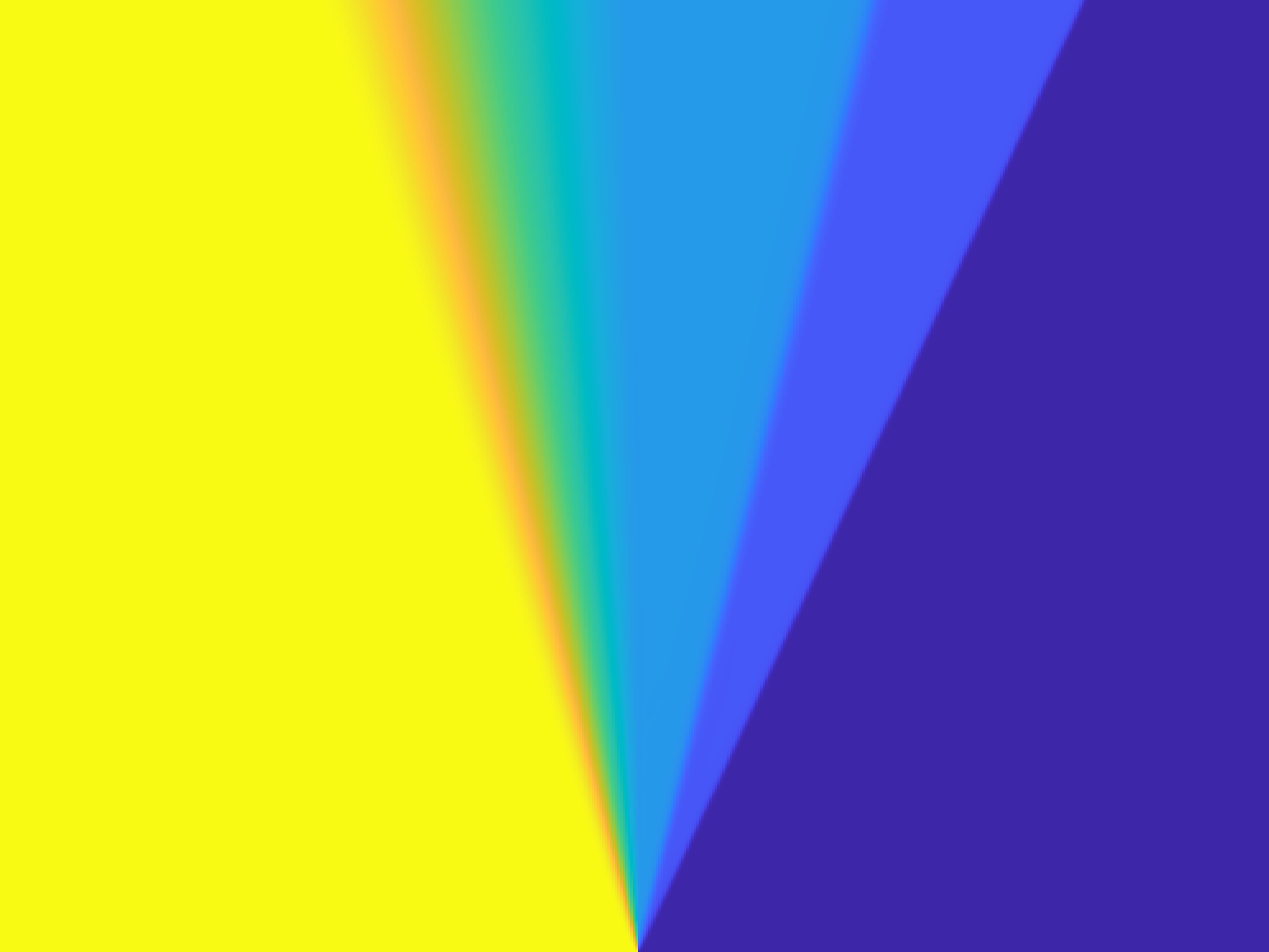}}
	~
	\subfloat[]{\includegraphics[width=.49\textwidth,trim={0mm 0mm 0mm 0mm},clip]{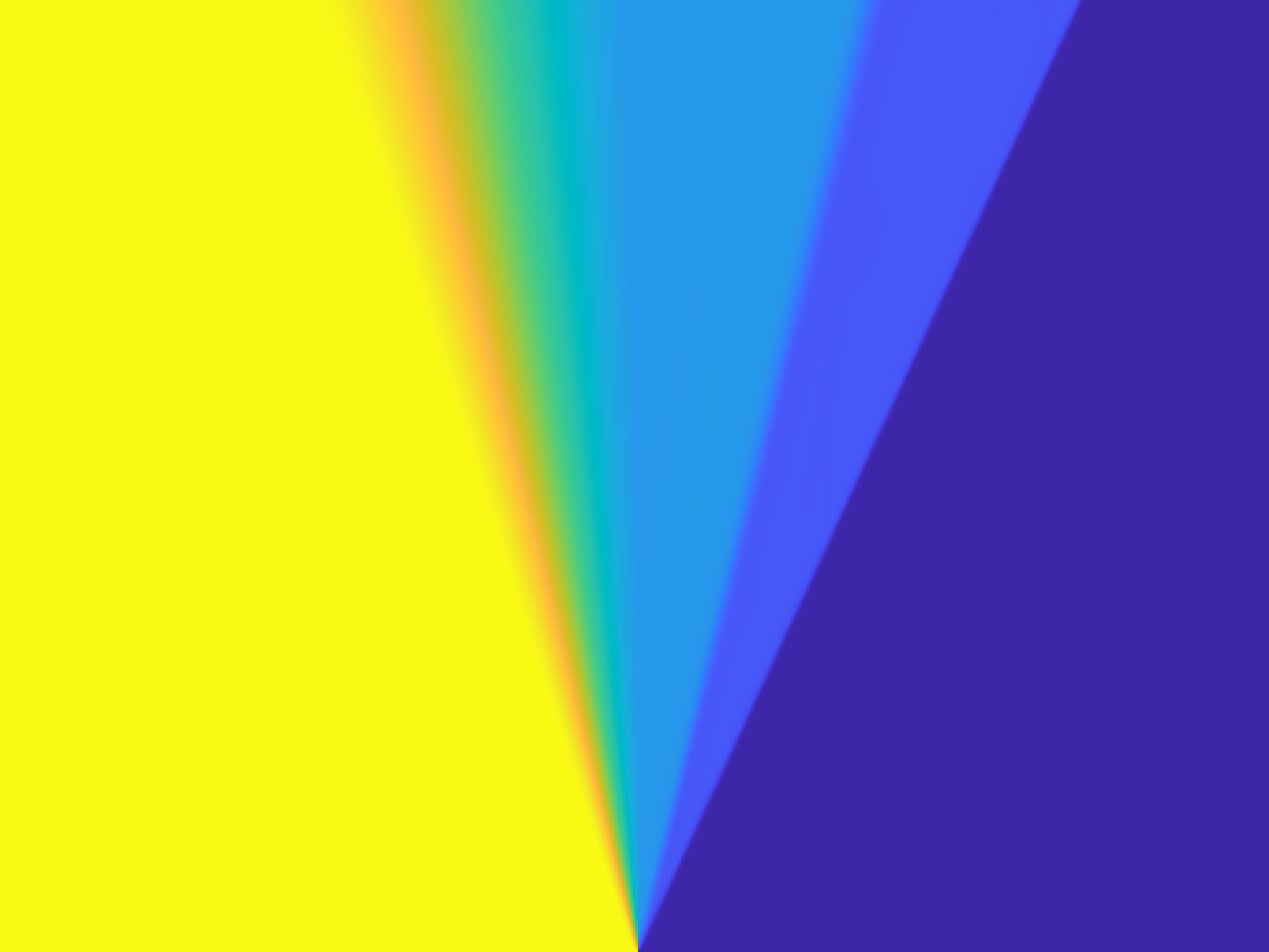}}
	\\
	\subfloat[]{\includegraphics[width=.49\textwidth,trim={0mm 0mm 0mm 0mm},clip]{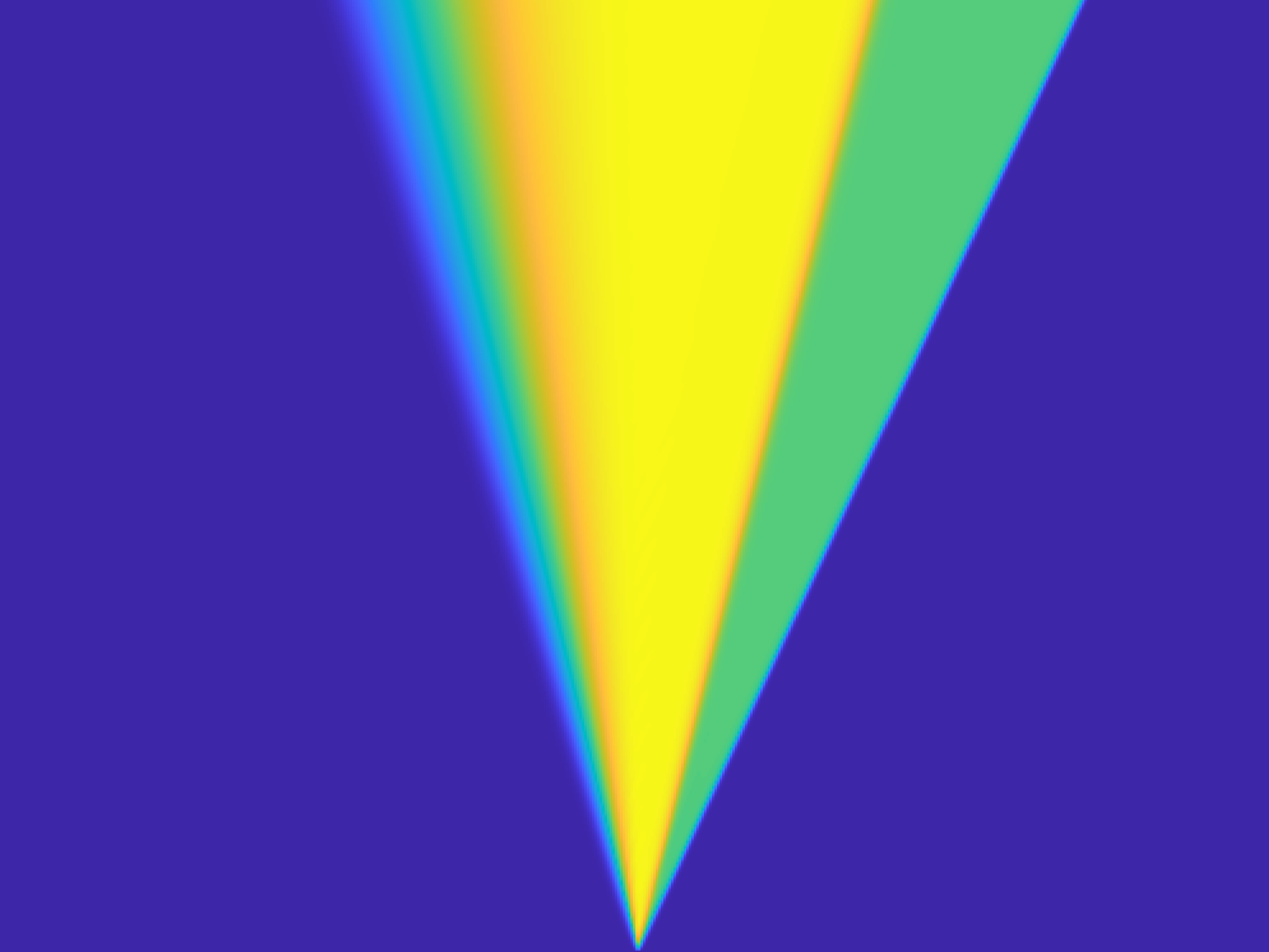}}
	~
	\subfloat[]{\includegraphics[width=.49\textwidth,trim={0mm 0mm 0mm 0mm},clip]{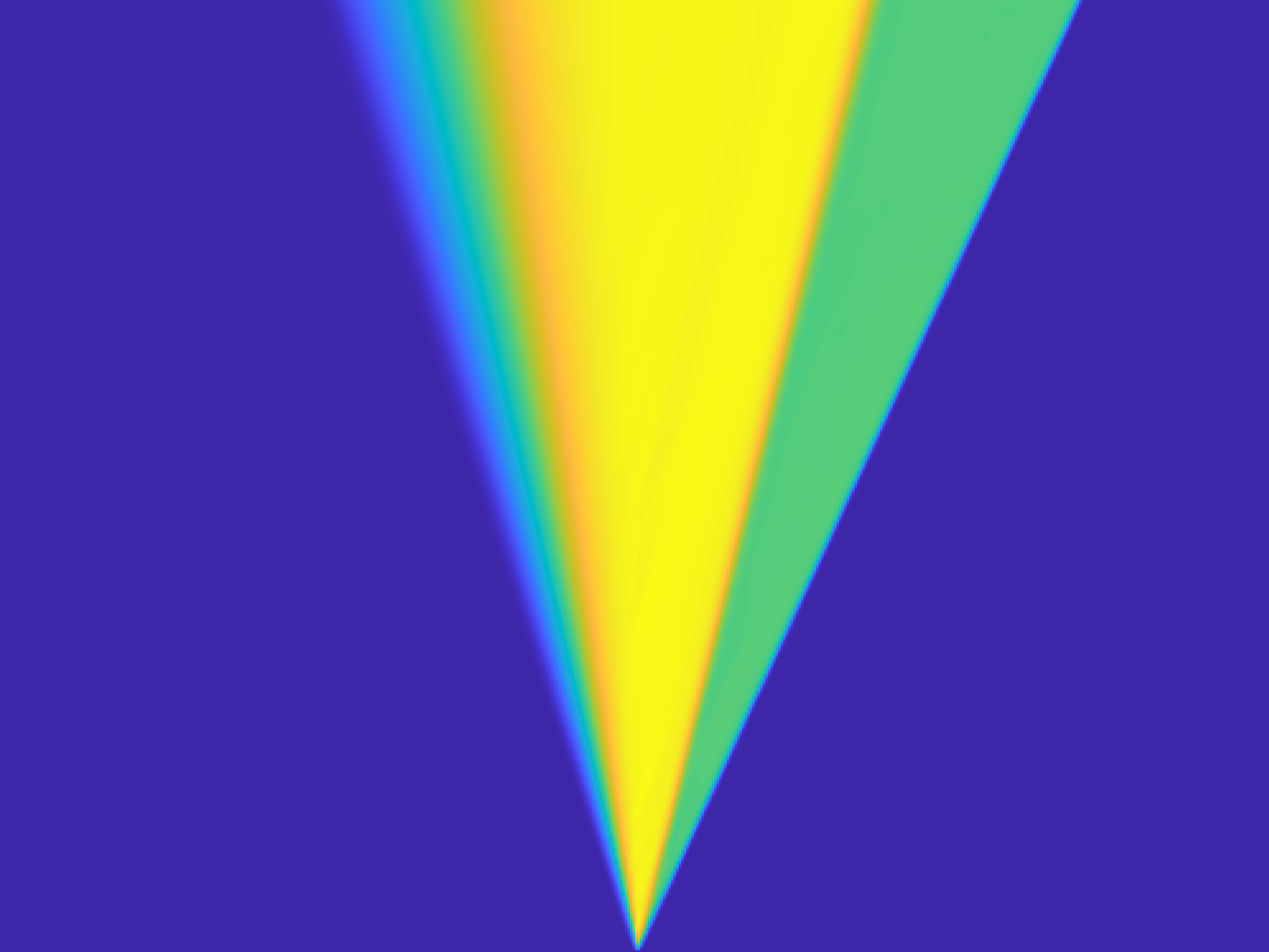}}
	\\
	\subfloat[]{\includegraphics[width=.49\textwidth,trim={0mm 0mm 0mm 0mm},clip]{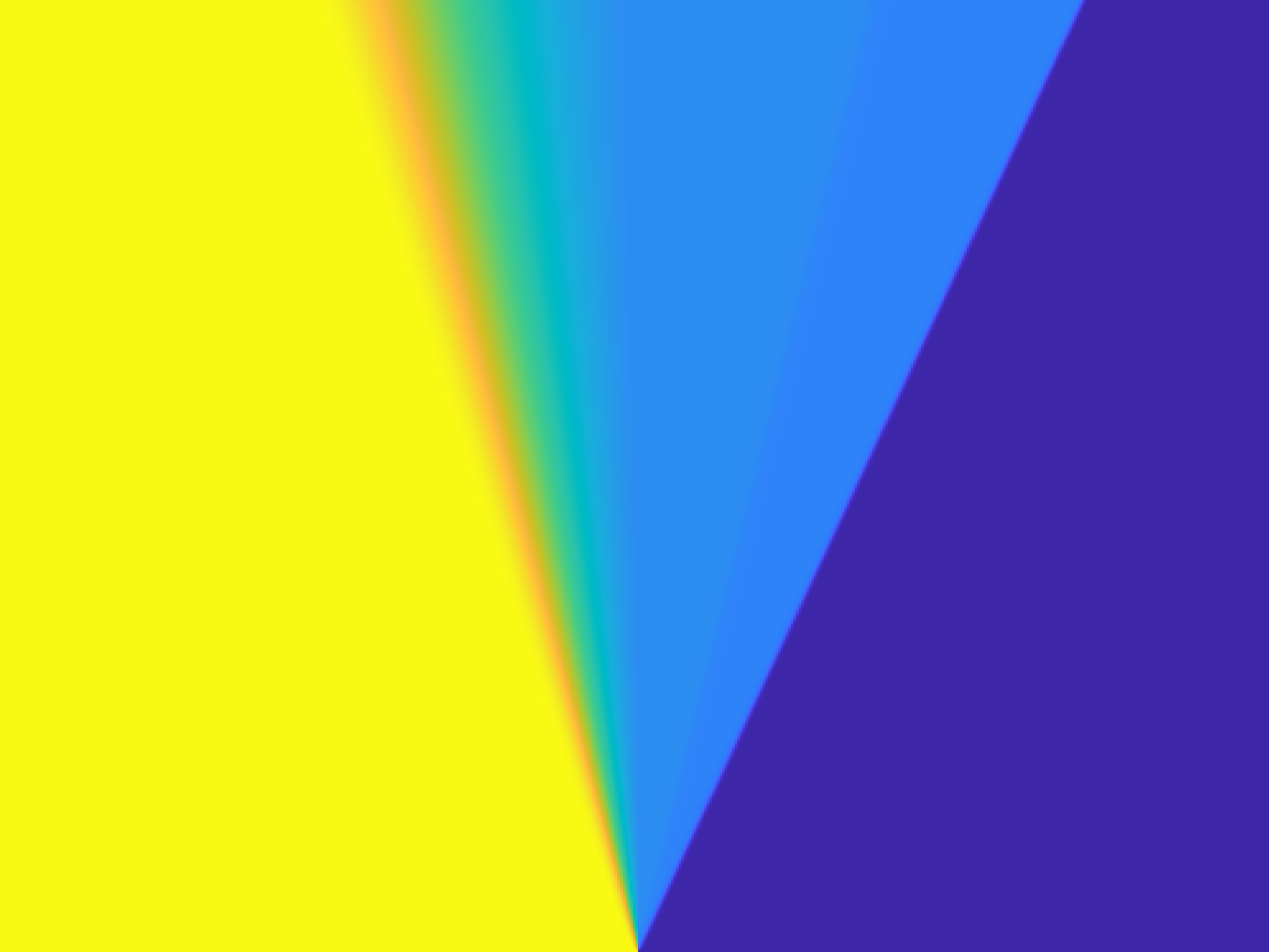}}
	~
	\subfloat[]{\includegraphics[width=.49\textwidth,trim={0mm 0mm 0mm 0mm},clip]{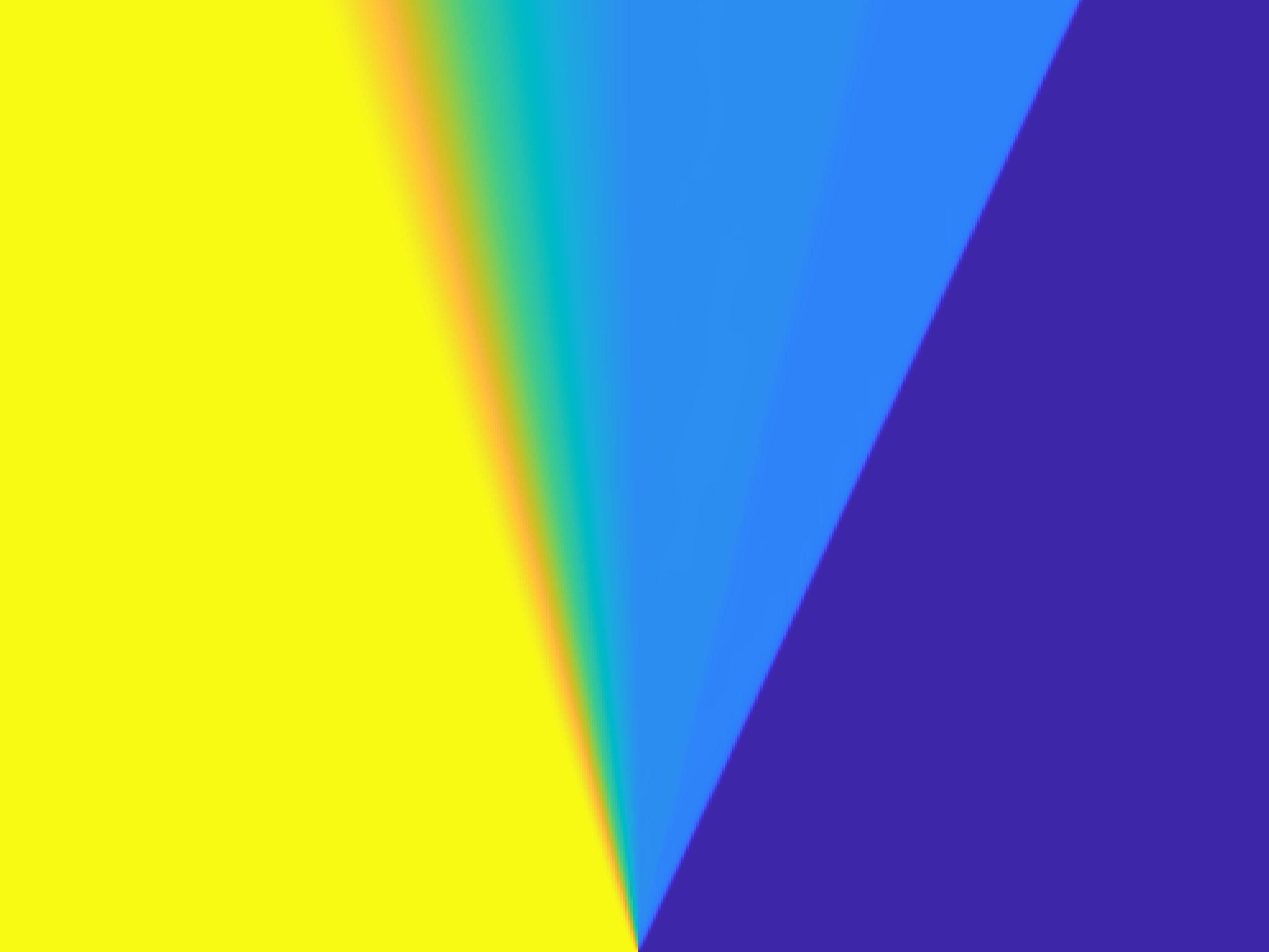}}
	\caption{Space-time snapshots of density (top), momentum (center) and energy (bottom) for a simulation only relying on full HDM solutions (left) and our AROM (right).}
	\label{fig:t14_snapshots_spaceTime_}
\end{figure}

Fig. \ref{fig:t14_map_} shows the points selected by sampling matrix $\hat{S}_k$. The first $w = 5$ snapshots are obtained using full HDM solves and, thus, are fully highlighted in yellow. From this figure, it can be noticed that the points are mainly concentrated on the propagating expansion, contact and shock waves. Sampling also takes place outside the range of influence of point $x = 0.5$. We can attribute this to the development of unfiltered non-physical structures.
From Fig. \ref{fig:t14_snapshot_} we can observe that a lower full HDM solve frequency leads the shock to lag behind.
The underestimation of the shock waves velocities can be attributed, at least in part, to the sampling algorithm relying on the solution on time instance $k$ to determine the sampling points at time $k+1$. This systematically leads the method to fail to sample regions that are dynamically relevant in the immediate future (e.g., downwind of the shock). A sampling strategy that aims at fixing this problem has been recently proposed \cite{singh2023lookahead}.
Figure \ref{fig:t14_coarse} provides a comparison of the AROM ($z = \infty$) and full HDM solutions on coarser grids. For this problem, the AROM is generally more accurate despite slightly underestimating the shock velocities.
Figure \ref{fig:t14_snapshots_spaceTime_} compares solutions between this AROM and a simulation relying only on full HDM solves. The AROM recovers the main features of the flow with small discrepancies in the range of influence of point $x = 0.5$.

\begin{figure}[hbt!]
	\centering
	\begin{tikzpicture}
\begin{groupplot} [
group style={group size = 2 by 1, horizontal sep = 1.8cm, vertical sep = 1.4cm}]
\nextgroupplot[width=.48\textwidth, xtick={3,4,5,6,7,8,9,10}, ytick={4,6,8,10,12,14,16}, xlabel={$m$}, ymax=16, xmax=10, ylabel={$\bar{s} (\%)$}, xmin=3, ymin=4]
\addplot [green, thick, mark options={solid, thin}, mark=o, mark size=3, mark repeat=0]
coordinates {
( 3.00000000e+00,  6.23451311e+00)
( 4.00000000e+00,  7.36543227e+00)};\label{line:fig_sod_wm_4}

\addplot [blue, thick, mark options={solid, thin}, mark=diamond, mark size=3, mark repeat=0]
coordinates {
( 3.00000000e+00,  5.86604070e+00)
( 4.00000000e+00,  8.11261802e+00)
( 5.00000000e+00,  8.85197248e+00)
( 6.00000000e+00,  9.46341581e+00)};\label{line:fig_sod_wm_6}

\addplot [red, thick, mark options={solid, thin}, mark=triangle, mark size=3, mark repeat=0]
coordinates {
( 3.00000000e+00,  5.98913257e+00)
( 4.00000000e+00,  9.37988201e+00)
( 5.00000000e+00,  9.80076385e+00)
( 6.00000000e+00,  1.09356991e+01)
( 7.00000000e+00,  1.15194316e+01)
( 8.00000000e+00,  1.22481436e+01)};\label{line:fig_sod_wm_8}

\addplot [yellow, thick, mark options={solid, thin}, mark=triangle, mark size=3, mark repeat=0]
coordinates {
( 3.00000000e+00,  5.06403589e+00)
( 4.00000000e+00,  9.24956125e+00)
( 5.00000000e+00,  1.05744555e+01)
( 6.00000000e+00,  1.06031703e+01)
( 7.00000000e+00,  1.20025622e+01)
( 8.00000000e+00,  1.25804716e+01)
( 9.00000000e+00,  1.34700262e+01)
( 1.00000000e+01,  1.46031542e+01)};\label{line:fig_sod_wm_10}

\nextgroupplot[width=.48\textwidth, xtick={3,4,5,6,7,8,9,10}, ytick={.1,.2,.3,.4,.5,.6,.8,1}, xlabel={$m$}, ymax=0.6, xmax=10, ylabel={$\bar{e} (\%)$}, xmin=3, ymin=0.1]
\addplot [green, thick, mark options={solid, thin}, mark=o, mark size=3, mark repeat=0]
coordinates {
( 3.00000000e+00,  4.03040035e-01)
( 4.00000000e+00,  4.12299744e-01)};\label{line:fig_sod_wm_4}

\addplot [blue, thick, mark options={solid, thin}, mark=diamond, mark size=3, mark repeat=0]
coordinates {
( 3.00000000e+00,  3.66201022e-01)
( 4.00000000e+00,  3.46259485e-01)
( 5.00000000e+00,  2.93301374e-01)
( 6.00000000e+00,  2.68118224e-01)};\label{line:fig_sod_wm_6}

\addplot [red, thick, mark options={solid, thin}, mark=triangle, mark size=3, mark repeat=0]
coordinates {
( 3.00000000e+00,  4.57329738e-01)
( 4.00000000e+00,  3.91168393e-01)
( 5.00000000e+00,  2.79701164e-01)
( 6.00000000e+00,  2.30623136e-01)
( 7.00000000e+00,  2.03686085e-01)
( 8.00000000e+00,  1.92264300e-01)};\label{line:fig_sod_wm_8}

\addplot [yellow, thick, mark options={solid, thin}, mark=triangle, mark size=3, mark repeat=0]
coordinates {
( 3.00000000e+00,  5.76528597e-01)
( 4.00000000e+00,  4.93833229e-01)
( 5.00000000e+00,  3.39600736e-01)
( 6.00000000e+00,  2.79896589e-01)
( 7.00000000e+00,  2.17549264e-01)
( 8.00000000e+00,  1.83409805e-01)
( 9.00000000e+00,  1.67141905e-01)
( 1.00000000e+01,  1.62786160e-01)};\label{line:fig_sod_wm_10}

\end{groupplot}\end{tikzpicture}
	\caption{Time averages of relative sampling (left) and relative error (right) as a function of the reduced-order dimension $m$ for windows of size $w = 4$ (\ref{line:fig_sod_wm_4}), $w = 6$ (\ref{line:fig_sod_wm_6}), $w = 8$ (\ref{line:fig_sod_wm_8}) and $w = 10$ (\ref{line:fig_sod_wm_10}).}
	\label{fig:t14_wm_study_}
\end{figure}
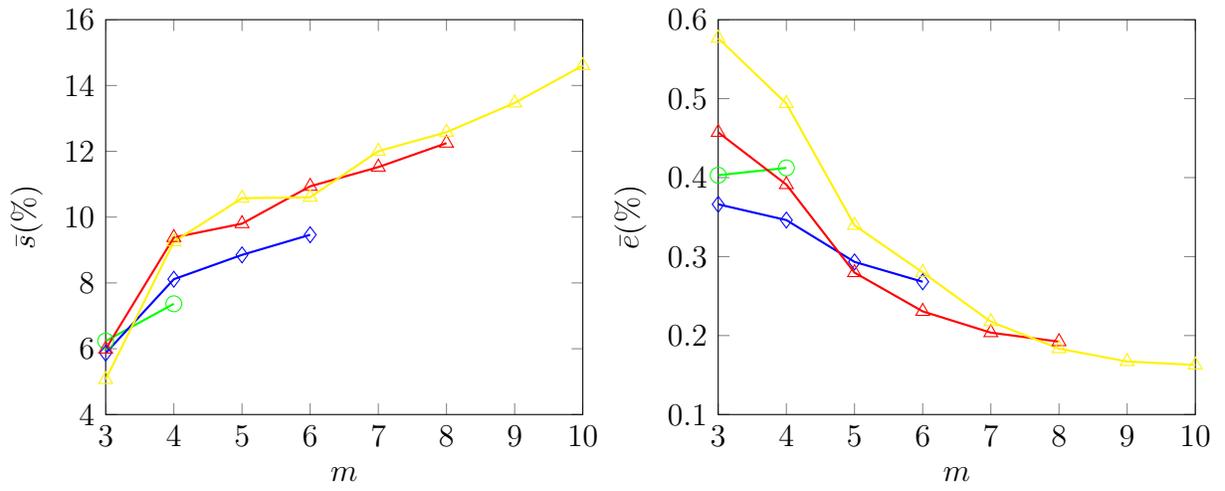

Figure \ref{fig:t14_wm_study_} shows time average error and sampling responses to different values of window size $w$ and number of POD modes $m$.
For $w = 4$, the error is the smallest for $m = 3$ despite the smaller sampling. An additional mode degrades the solution by adding non-physical structures that lead to bigger sampling matrices that are not completely dissipated by the filters.
The opposite trend is observed for all other cases. A bigger basis generally leads to a more accurate solutions at the cost of bigger sampling matrices. However, these ROMs are considerably less accurate if the number of modes used in the reconstruction is too small. 
Furthermore, we can observe that larger windows and bigger bases lead to bigger sampling matrices. This is expected as larger basis results in additional ODEIM points and is less of an issue for multidimensional problems because they usually lead to sparser sampling  (Section~\ref{sec:implosion}). 
As discuss in Section \ref{sec:complexity}, the cost of performing POD is also a quadratic function of window width $\mathcal{O} (N w^2)$. Therefore, a narrower window is preferred if the benefits of a larger window is negligible.

\begin{remark}
	Speedup evaluations for this problem are absent. First, high-dimensional models of one-dimensional problems can, for the most part, be easily solved using a laptop and, thus, reduced-order models are unnecessary. 
	Second, the timing of small problems are not representative of large-scale problems. For example, the cost of solving a system of nonlinear equations relative to the cost of other operations (e.g., residual and Jacobian evaluations) is typically much higher for large-scale problems.
	We only use this problem as proof of concept.
\end{remark}

\subsection{Model implosion}
\label{sec:implosion}

In this problem, we consider the two-dimensional ($d=2$) Euler equations in the domain $\Omega \subset (0,0.3)^2$ over the time interval $\mathcal{T} = (0, .5)$ with ratio of specific heats $\gamma = 1.4$ and initial condition (in terms of primitive variables) as

\begin{align}
	\rho (x,0) = 
	\begin{cases}
		\rho_\mathrm{in} & x \in \mathcal{D}\\
		\rho_\mathrm{out} & x \notin \mathcal{D}
	\end{cases}
	,
	\quad
	u(x,0) = (0,0),
	\quad
	P(x,0) = 
	\begin{cases}
		P_\mathrm{in} & x \in \mathcal{D}\\
		P_\mathrm{out} & x \notin \mathcal{D}
	\end{cases}
	\mbox{ .}
\end{align}
where $\rho_\mathrm{in} = 0.125$ and $P_\mathrm{in} = 0.14$ are the pressure and density inside the region $\mathcal{D} = \{ x \in \Omega \mid x_1 + x_2  \leq 0.15 \} \subset \Omega$ and $\rho_\mathrm{out} = 1$ and $P_\mathrm{out} = 1$ are the pressure and density outside $\mathcal{D}$. All four boundaries are taken to be walls, which causes the waves to reflect back into the domain when they reach a boundary. This is a model of an implosion that was adapted from \cite{Ghosh_2012_implosion}.

We solve this problem using a $100 \times 100$ uniform cartesian grid. We partition the time domain into $N_t = 1\mbox{,}650$ time steps (chosen for global accuracy and steep shock approximations).
As the previous problem, we filter hybrid solutions by sequentially applying second-,fourth- and sixth-order filters.
The full HDM frequency parameter is $z = 7$ and the reconstruction error threshold is set at $\delta = 0.90$.
The number of snapshots used in the reduced basis reconstruction is $w = 6$ and the number of POD modes used in the reconstruction is $m = 4$. 
For these parameters, $2 \leq J \leq 4$ with the average number of subiterations being $\bar{J} = 2.73$.

The time average sampling values are $\bar{p} = 0.23 \%$, $\bar{s} = 18.35 \%$ and $\bar{s}^{*} = 4.07\%$, and the hybrid snapshot sampling never exceeds $34\%$. Finally, we use speedup Eq. \eqref{eq:speedup} to evaluate the relative execution time performance of our AROM method. For this problem in particular, the speedup is $\mathcal{S} = 4.52$.

\begin{figure}[hbt!]
	\centering
	\subfloat[]{\includegraphics[width=.27\textwidth,trim={0mm 0mm 0mm 0mm},clip]{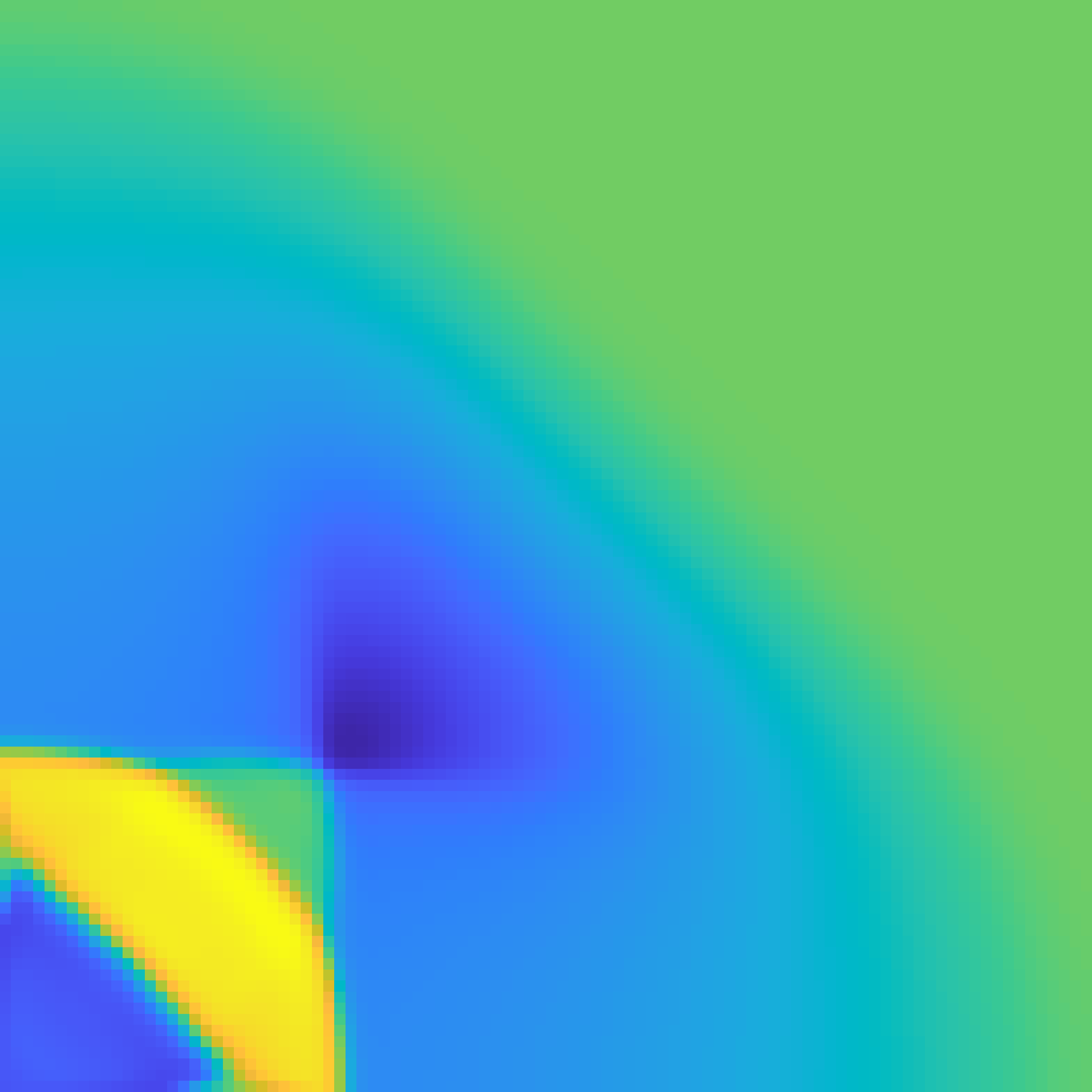}}
	~
	\subfloat[]{\includegraphics[width=.27\textwidth,trim={0mm 0mm 0mm 0mm},clip]{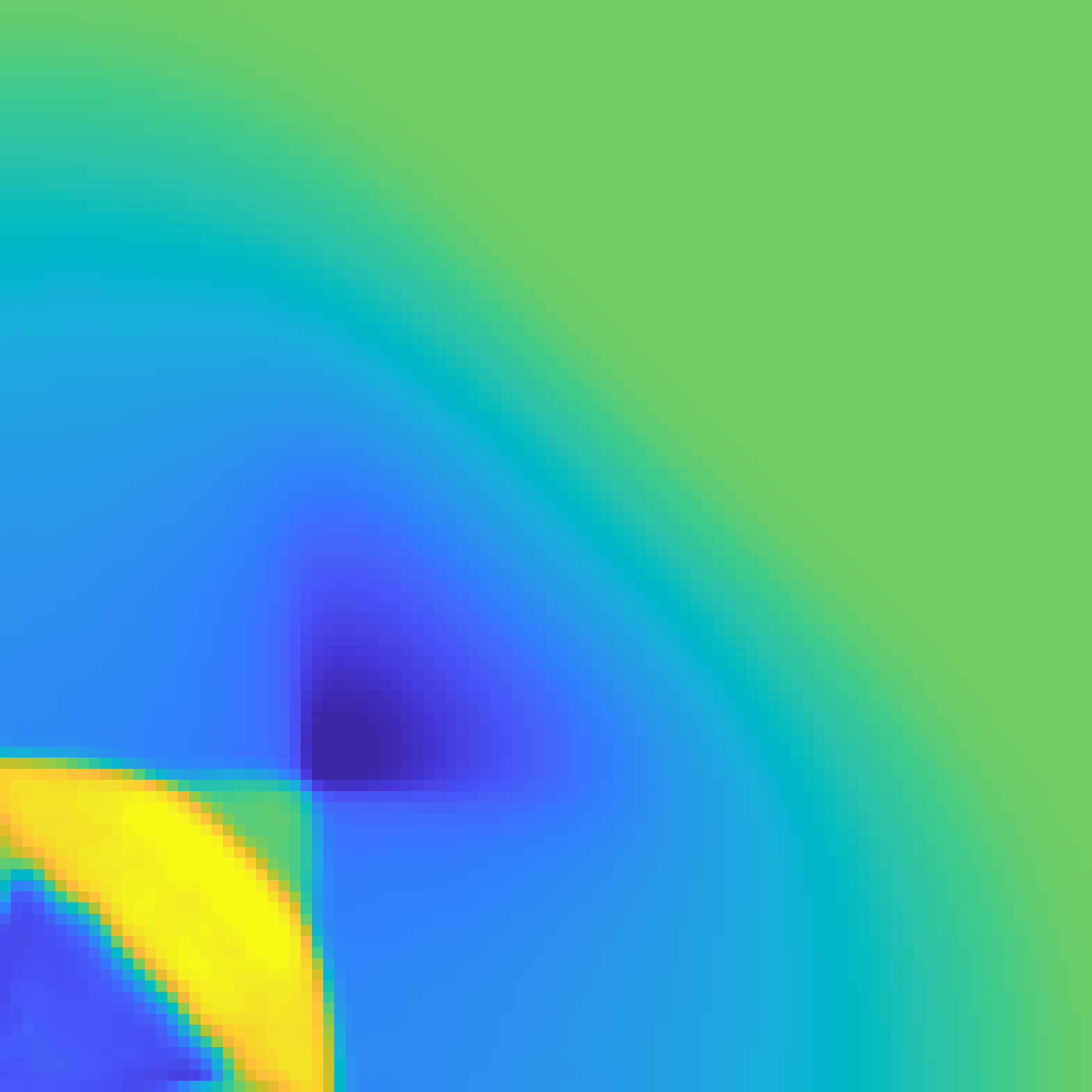}}
	~
	\subfloat[]{\includegraphics[width=.27\textwidth,trim={0mm 0mm 0mm 0mm},clip]{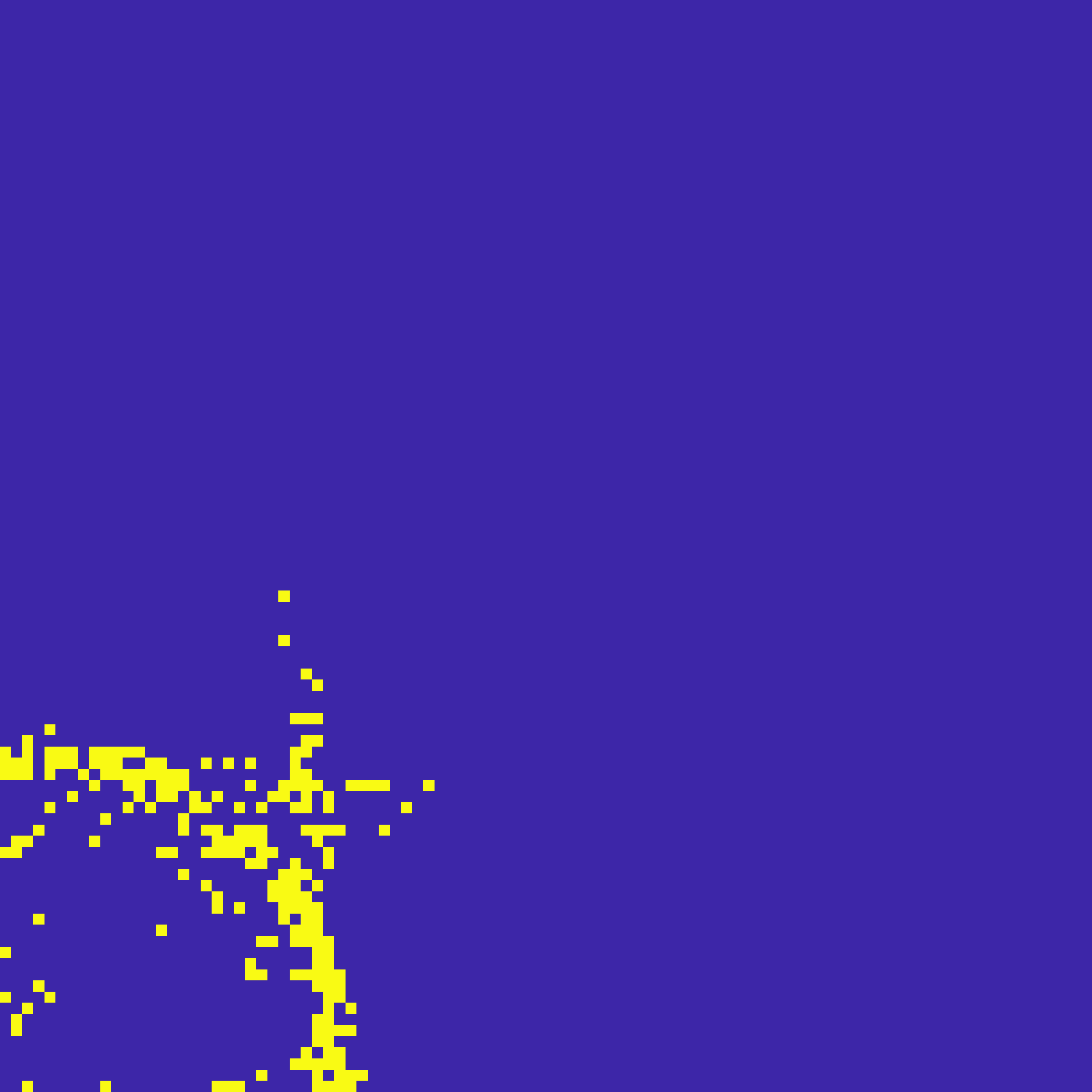}}
	\\
	\subfloat[]{\includegraphics[width=.27\textwidth,trim={0mm 0mm 0mm 0mm},clip]{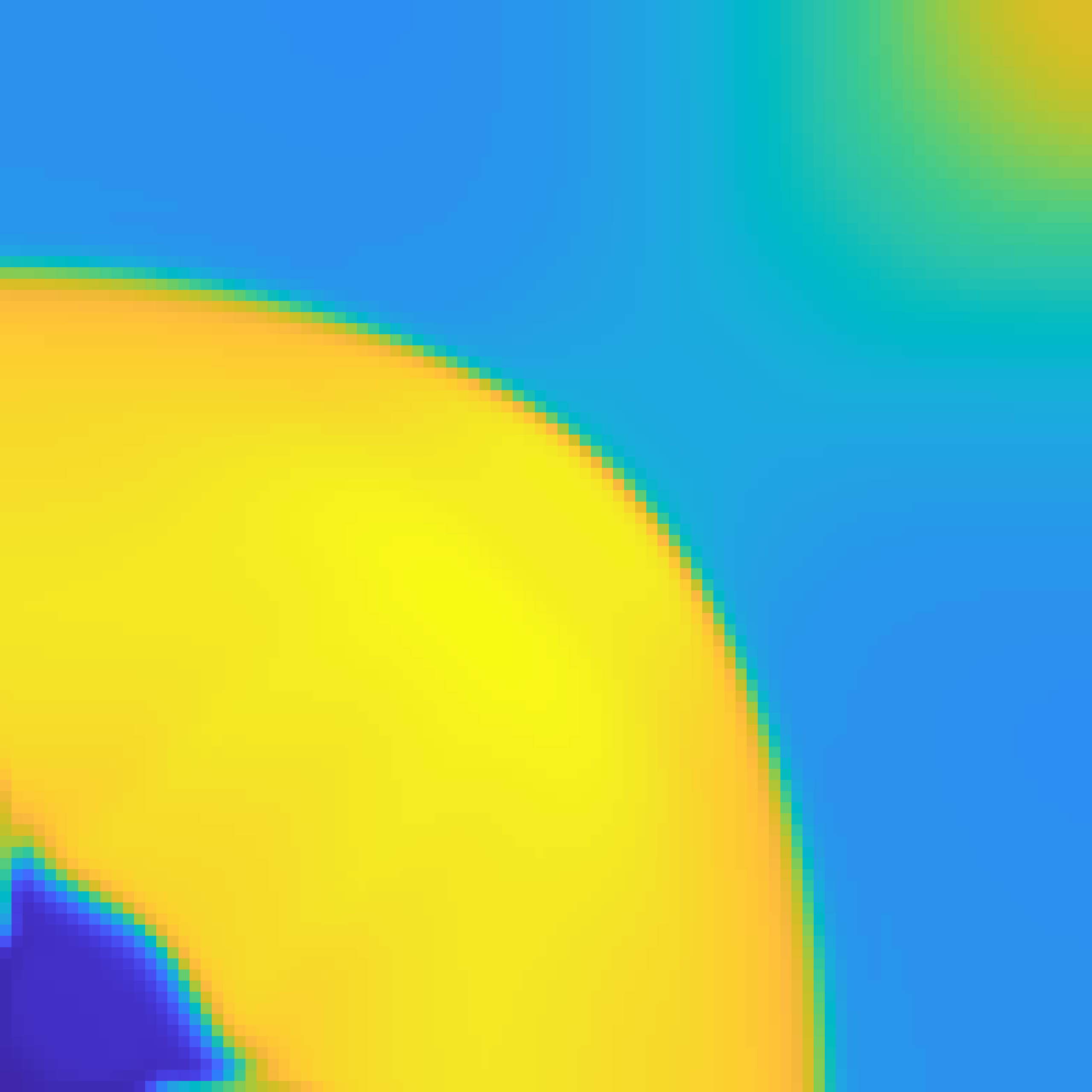}}
	~
	\subfloat[]{\includegraphics[width=.27\textwidth,trim={0mm 0mm 0mm 0mm},clip]{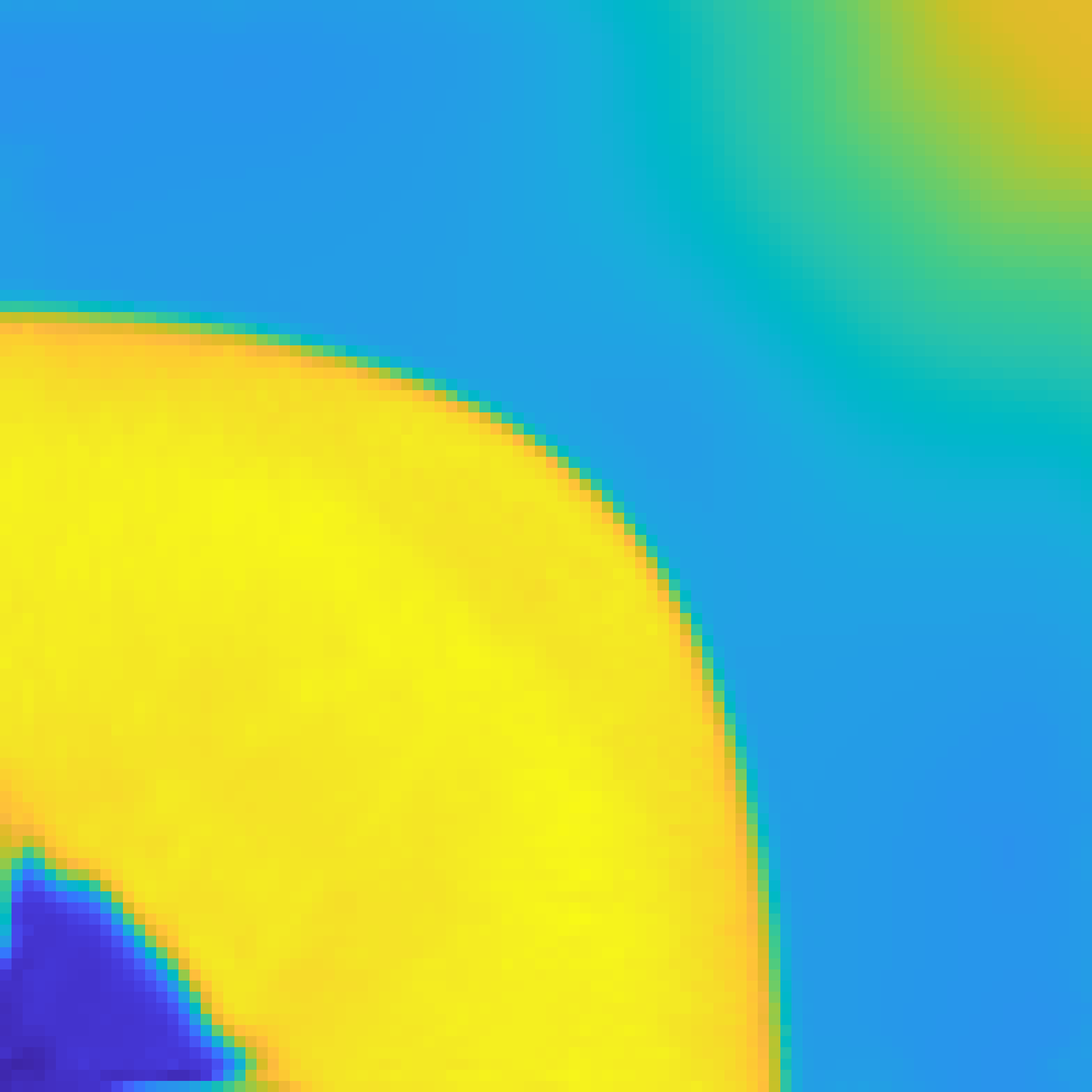}}
	~
	\subfloat[]{\includegraphics[width=.27\textwidth,trim={0mm 0mm 0mm 0mm},clip]{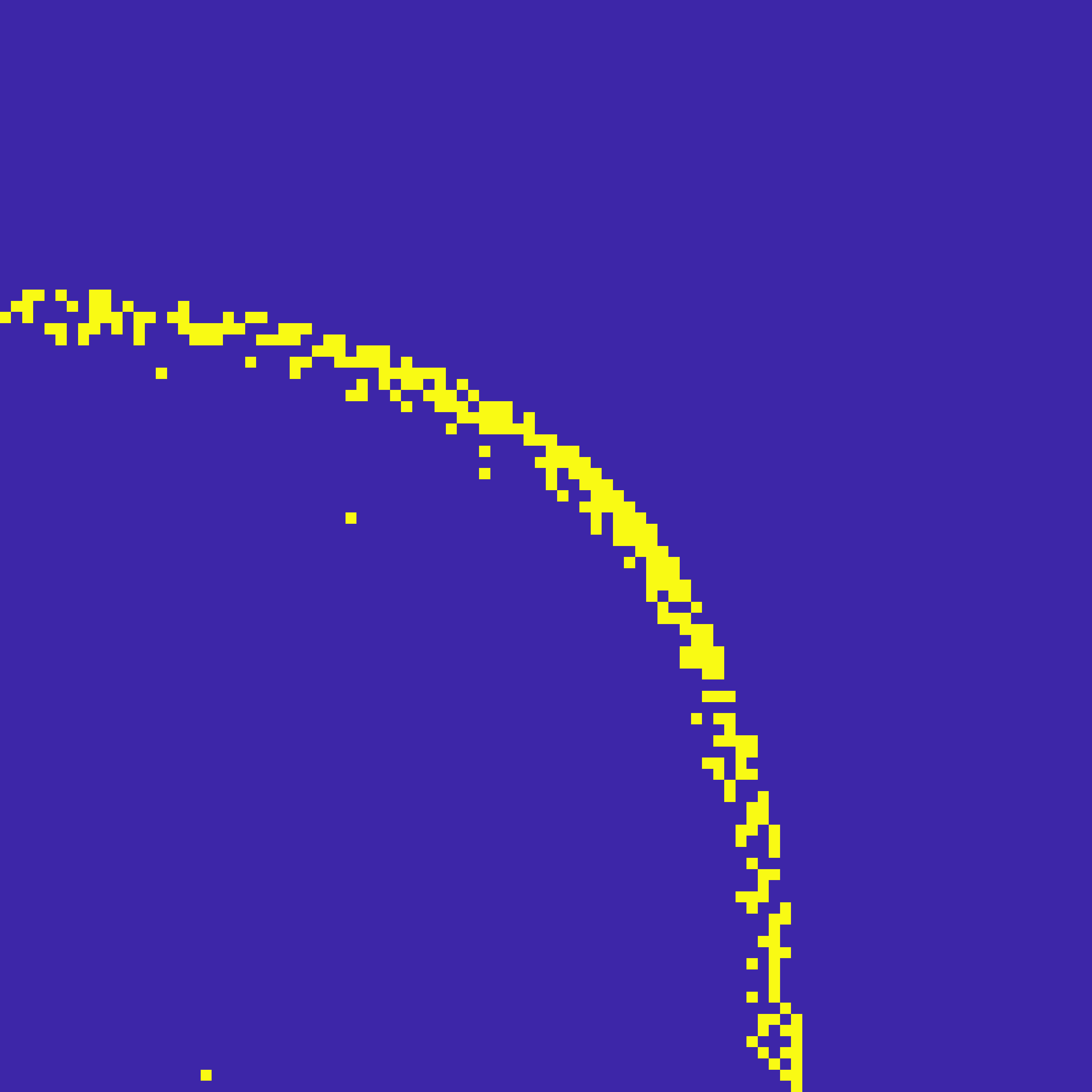}}
	\\
	\subfloat[]{\includegraphics[width=.27\textwidth,trim={0mm 0mm 0mm 0mm},clip]{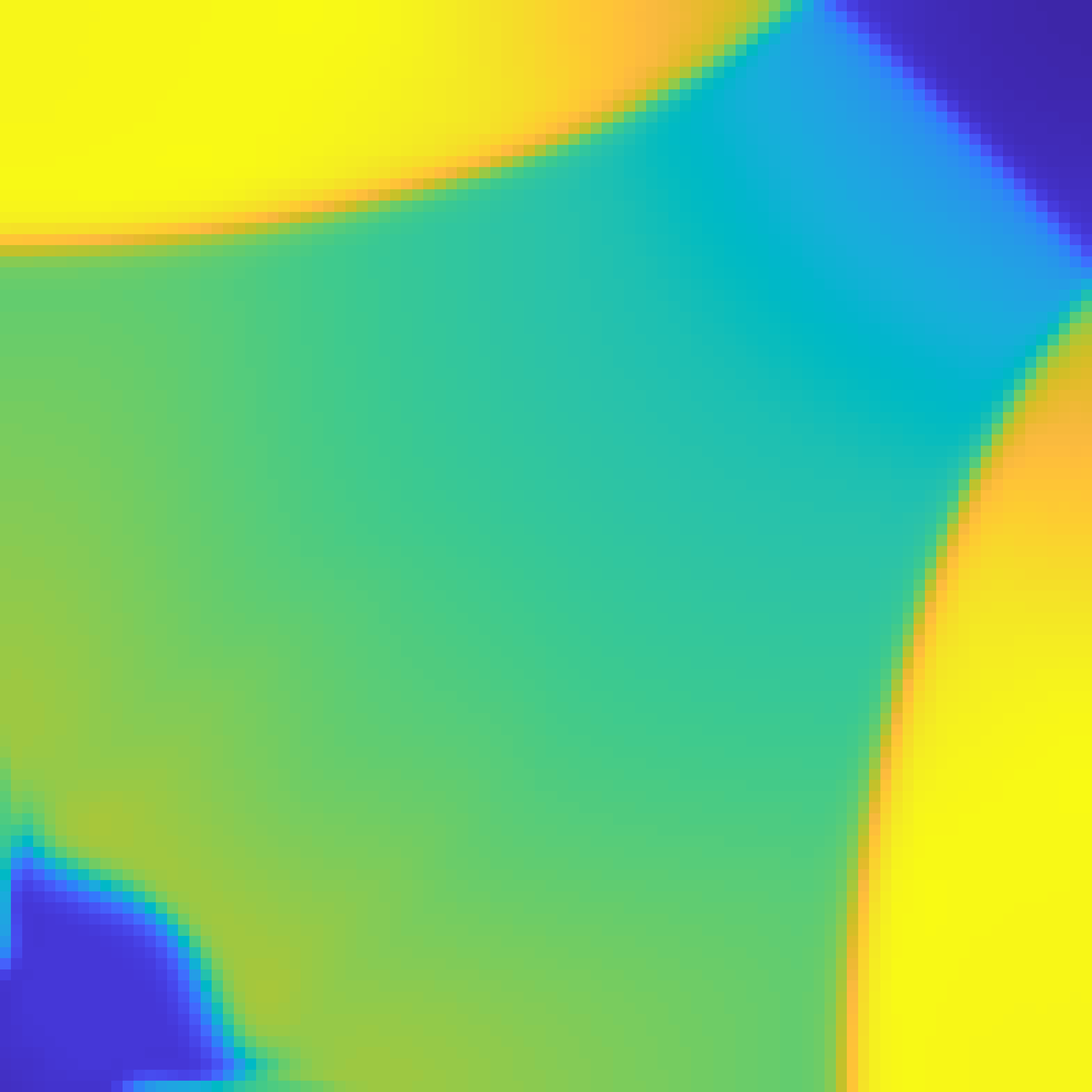}}
	~
	\subfloat[]{\includegraphics[width=.27\textwidth,trim={0mm 0mm 0mm 0mm},clip]{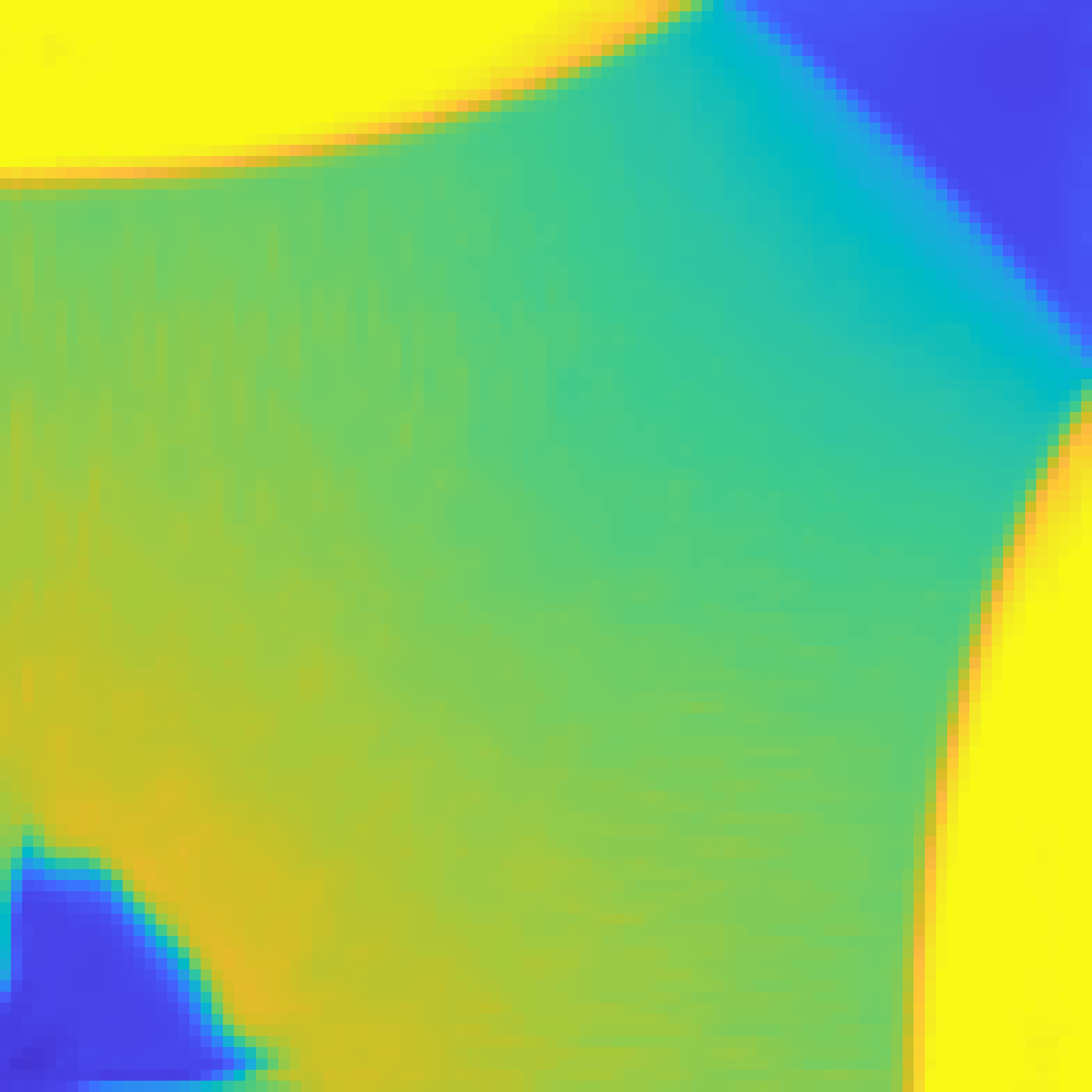}}
	~
	\subfloat[]{\includegraphics[width=.27\textwidth,trim={0mm 0mm 0mm 0mm},clip]{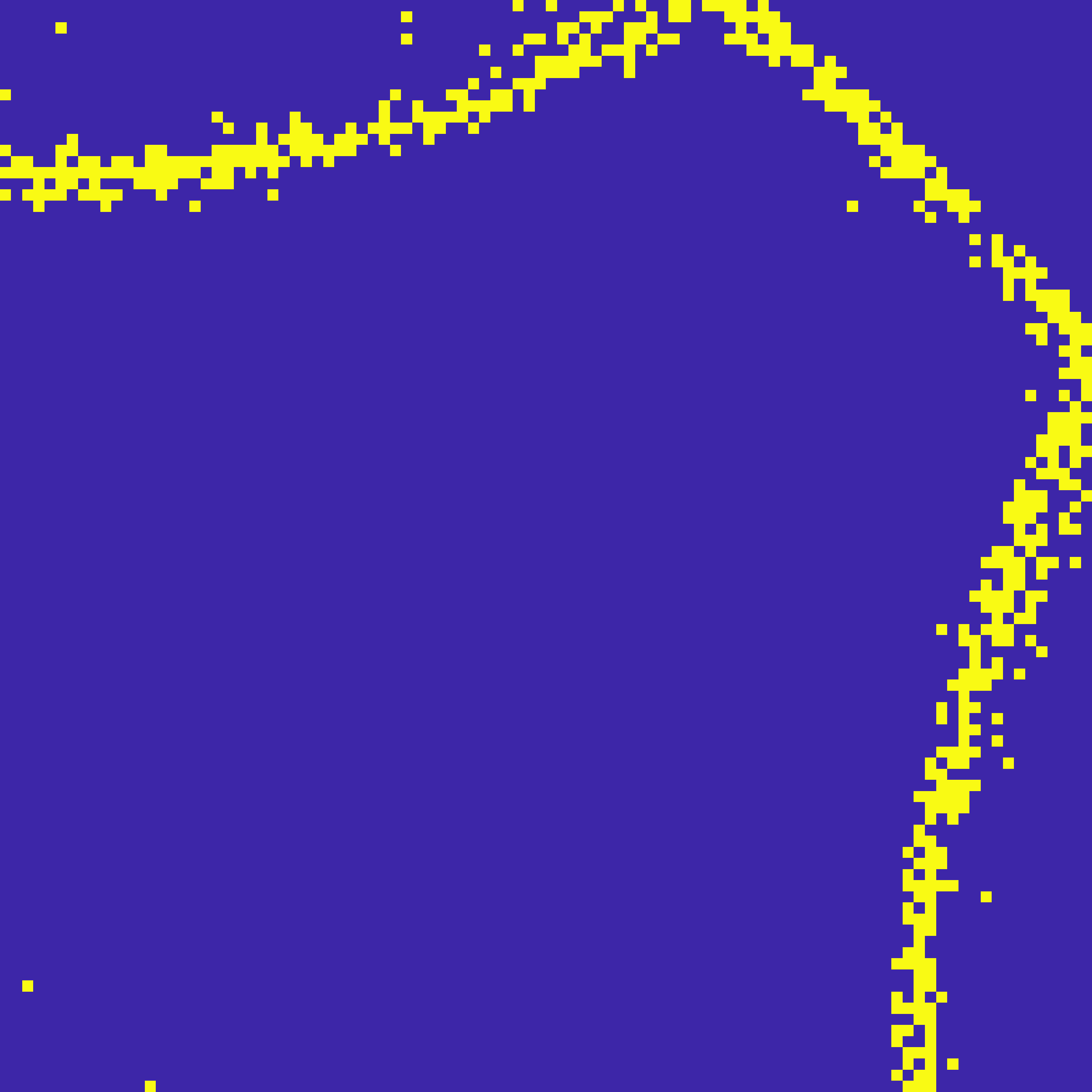}}
	\\
	\subfloat[]{\includegraphics[width=.27\textwidth,trim={0mm 0mm 0mm 0mm},clip]{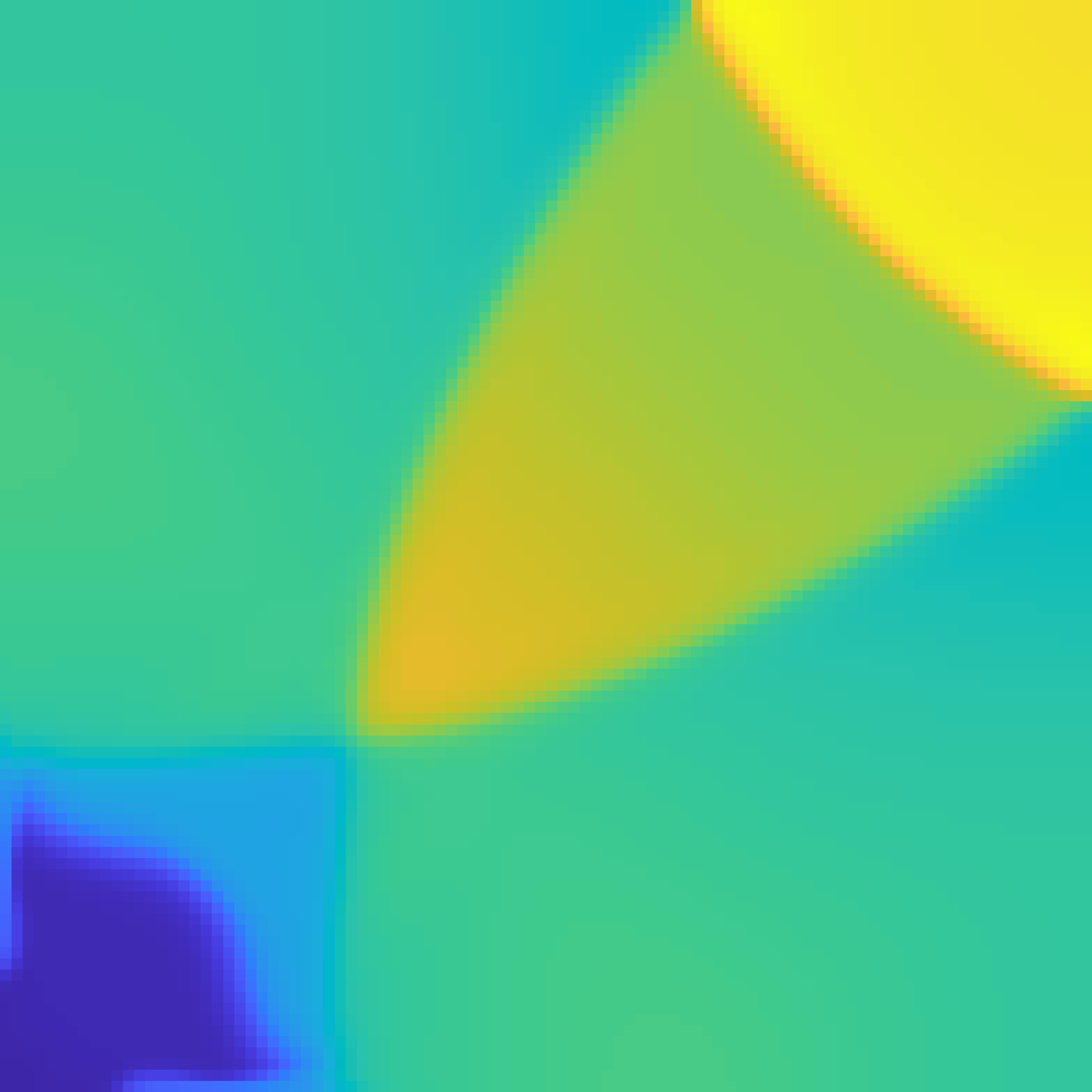}}
	~
	\subfloat[]{\includegraphics[width=.27\textwidth,trim={0mm 0mm 0mm 0mm},clip]{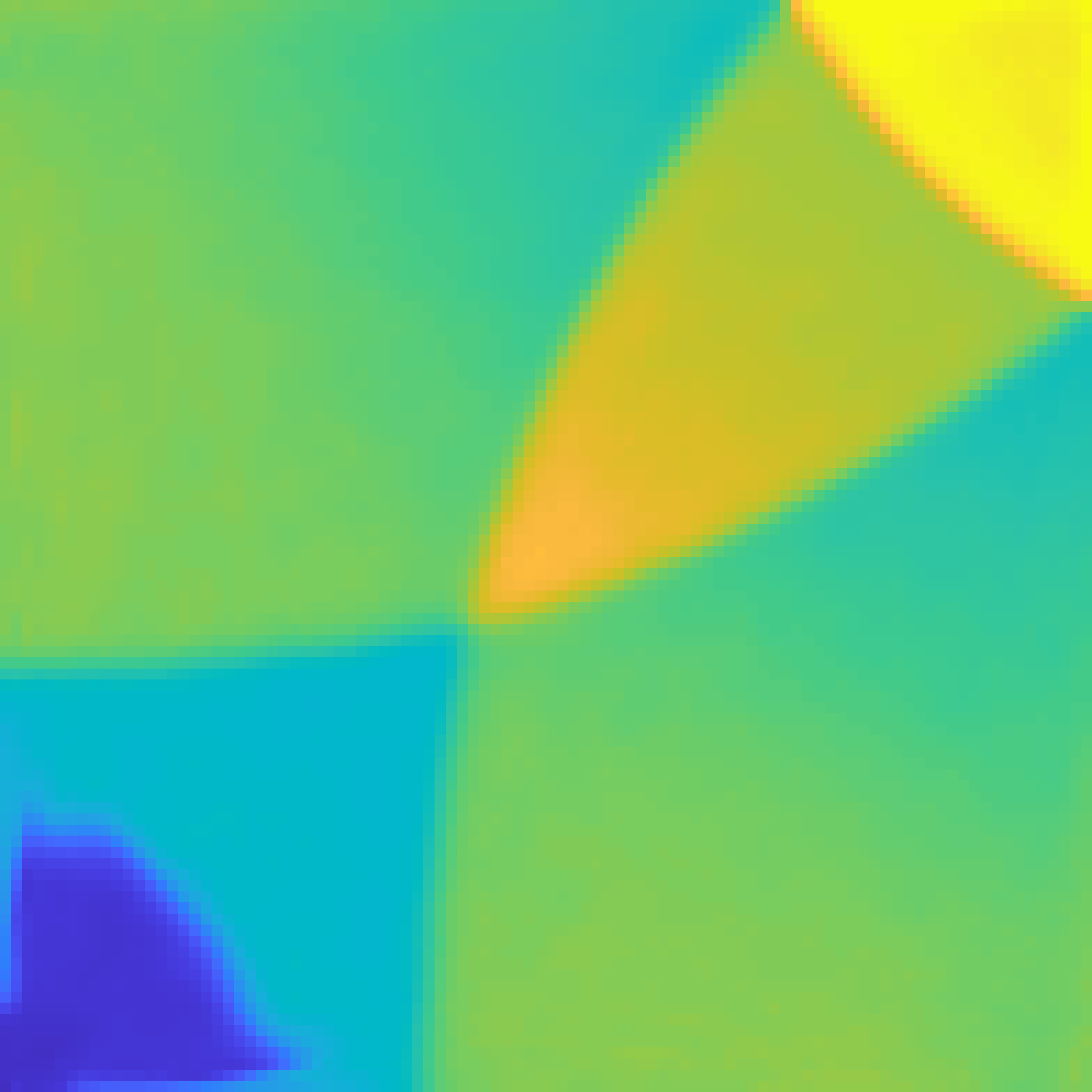}}
	~
	\subfloat[]{\includegraphics[width=.27\textwidth,trim={0mm 0mm 0mm 0mm},clip]{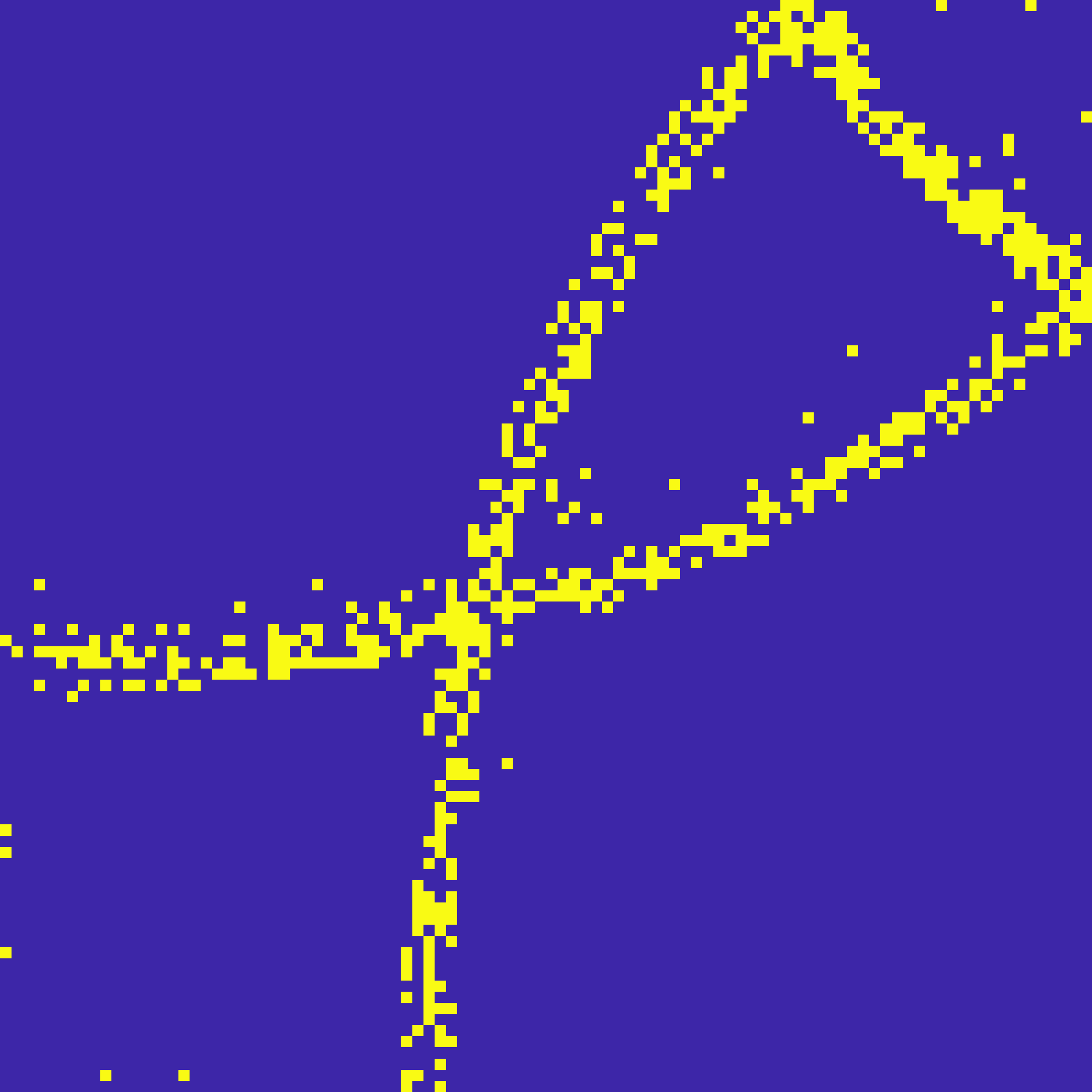}}
	\caption{Density snapshots a simulation only relying on full HDM solutions (left) and our AROM (center), and the sampling points corresponding to matrix $\hat{S}$ (right) at time instances $t = T/4$, $t = T/2$, $t = 3T/4$ and $t = T$ (top-to-bottom).}
	\label{fig:t19_snapshots_}
\end{figure}

\begin{figure}[hbt!]
	\centering
	\subfloat[]{\includegraphics[width=.33\textwidth,trim={0mm 0mm 0mm 0mm},clip]{T19_new/T19_q_hdm_4.png}}
	~
	\subfloat[]{\includegraphics[width=.33\textwidth,trim={0mm 0mm 0mm 0mm},clip]{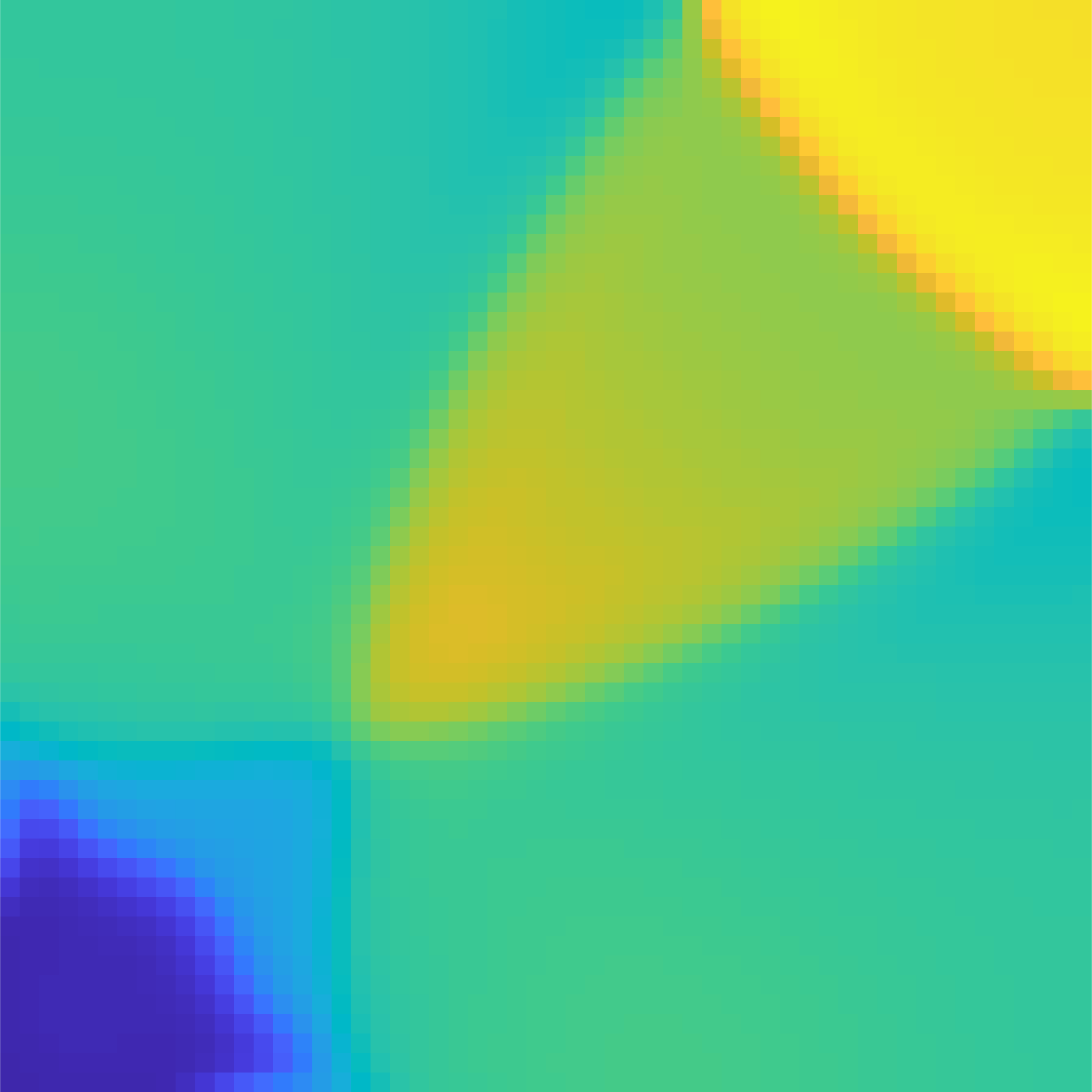}}
	~
	\subfloat[]{\includegraphics[width=.33\textwidth,trim={0mm 0mm 0mm 0mm},clip]{T19_new/T19_q_arom_4.png}}
	\caption{Density snapshots a simulation only relying on full HDM solutions (left and center) and our AROM (right) at time instances $t = T$. The center snapshot uses a coarser grid ($58 \times 58$).}
	\label{fig:t19_snaepshots_coarse}
\end{figure}

Figure \ref{fig:t19_snapshots_} shows snapshots of a simulation relying only on full HDM solves, our AROM, and the cells selected by sampling matrix $\hat{S}_k$.
For all four time instances, the AROM is capable of solving the main features of the problem with only some minor discrepancies.
Larger errors are observed near boundaries, which can be addressed by separately sampling the boundaries and interior.
We compare in Fig. \ref{fig:t19_snaepshots_coarse} the AROM to a simulation relying only on full HDM solves on a coarser grid of equivalent cost.
Similar to the previous problem, the AROM solution is lagged relative to the fine-grid HDM solution. Still the AROM solution diffusion error is considerably smaller than the coarse-grid HDM solution where shocks are blurred and features underresolved.
%Similarly the previous problem, the AROM solution is also lagged but to a greater degree. In a coarser grid, failing to sample a single element can have a bigger impact in the immediate future. Moreover, we can also expect waves traveling at higher speeds to be more affected.
%
%On the positive side, the AROM solution diffusion error is also considerably smaller.

\section{Conclusions and future directions}
\label{sec:conclusions}

In this work, an adaptive reduced-order model is applied to convection-dominated problems. 
This approach relies on local HDM solves to obtain an accurate representation of the main flow features. The remainder of the flow is represented using a subspace approximation trained using previous snapshots.
The performance of the our approach is validated on two compressible flow problems with moving sharp gradient features. The first is the one-dimensional canonical Sod's shock tube problem, which it is used to conduct a parametric study. The second is a considerably more challenging two-dimensional problem simulating an implosion inside a box.
Results show that the proposed method is capable of accelerating convection-dominated unsteady CFD problems. If the sampling matrices remain sufficiently small throughout the simulation, a brief complexity analysis establishes that the speedup depends mainly on the full HDM solution frequency parameter $z$.  
Our first test case demonstrates that filtering combined with a residual error indicator allows for higher $z$ and, thus, is a crucial ingredient for cheaper and accurate AROMs. 
%
%One contribution of this work and an important component of the proposed method is the dynamic sampling matrix $\hat{S}_k$. Our time adaptive approach selects the smallest sampling set satisfying a predefined error tolerance at each time instance. This allows the sampling matrix to shrink or expand in an attempt to avoid undersampling and oversampling.
% 
Furthermore, the shock tube problem shows that narrower windows and smaller bases are sufficient to generate cheap and accurate AROMs.
% 
%Another important contribution is the partial HDM sampling used to construct the hybrid snapshots. It requires an approximation to the state on cells neighboring sample points, which can be made more accurate through subiterations and generally results in some accuracy gain without a significant cost increase.

The method could benefit from further research in multiple ways.
First, our current dynamic sampling procedure selects entries based only on their relative contribution to the total reconstruction error. For example, if the error tolerance is chosen to be too strict, this can lead to bigger sampling matrices than necessary if the residual is uniformly distributed across the mesh. Therefore, a better sampling algorithm could improve robustness and decrease cost.
Another research direction is boundary sampling. As previously discussed, accuracy at the boundaries could possibly be improved with little effort by sampling interior and boundary cells separately.
Finally, our approach relies on linear order reduction for most hybrid snapshots entries, i.e., the adapted basis $\Phi_k$ is used to compute the solution at the $\breve{S}_k$ indices. We avoid the Kolmogorov $n$-width problem by relying on the local low-rank structure of convection-dominated problems. Unfortunately, ROMs built on POD can struggle in predictive settings for even very simple problems. Nonlinear model reduction techniques could potentially overcome this barrier and produce AROMs less dependent on full HDM solves.

\section*{Data sharing}
Data sharing is not applicable to this article as no new data were created or analyzed in this study.

\section*{Acknowledgments}

This material is based upon work supported by the Air Force Office of Scientific Research (AFOSR) under award numbers FA9550-20-1-0236 and FA9550-22-1-0004. The content of this publication does not necessarily reflect the position or policy of any of these supporters, and no official endorsement should be inferred.

%% If you have bibdatabase file and want bibtex to generate the
%% bibitems, please use
%%
 \bibliographystyle{elsarticle-num} 
 \bibliography{cas-refs}

%% else use the following coding to input the bibitems directly in the
%% TeX file.

% \begin{thebibliography}{00}

% %% \bibitem{label}
% %% Text of bibliographic item

% \bibitem{}

% \end{thebibliography}
\end{document}